 \documentclass[psamsfonts]{article}
 
\usepackage{amsfonts}
\usepackage{amsmath}
\usepackage{amssymb}
\usepackage{amscd}
\usepackage{amstext}
\usepackage{amsthm}

 \title{ Full extensions and approximate unitary equivalences
 \thanks{Research partially supported by NSF grants DMS .
         AMS 2000 Subject Classification Numbers:
                         Primary 46L05,
                                46L35.
                        Key words:
                         Extension of $C^*$-algebras, simple $C^*$-algebras
                                                      \protect\\}}

\author{Huaxin Lin\\
Department of Mathematics\\
University of Oregon\\
Eugene, Oregon 97403-1222}
\begin{document}
\maketitle

\newcommand{\CA}{$C^*$-algebra}
\newcommand{\SCA}{$C^*$-subalgebra}
\newcommand{\aue}{approximate unitary equivalence}
\newcommand{\ayue}{approximately unitarily equivalent}
\newcommand{\mops}{mutually orthogonal projections}
\newcommand{\hm}{homomorphism}
\newcommand{\pisca}{purely infinite simple \CA}
\newcommand{\andeqn}{\,\,\,\,\,\, {\rm and} \,\,\,\,\,\,}
\newcommand{\QED}{\rule{1.5mm}{3mm}}
\newcommand{\morp}{contractive completely
positive linear
 map}
\newcommand{\asmorp}{asymptotic morphism}
\newcommand{\arrow}{\rightarrow}
\newcommand{\tdsum}{\widetilde{\oplus}}
\newcommand{\pa}{\|}  
\newcommand{\ep}{\varepsilon}
\newcommand{\id}{{\rm id}}
\newcommand{\aueeps}[1]{\stackrel{#1}{\sim}}
\newcommand{\aeps}[1]{\stackrel{#1}{\approx}}
\newcommand{\dt}{\delta}
\newcommand{\yu}{\fang}
\newcommand{\ca}{{\cal C}_1}
\newcommand{\Ad}{{\rm ad}}
\newtheorem{thm}{Theorem}[section]
\newtheorem{Lem}[thm]{Lemma}
\newtheorem{Prop}[thm]{Proposition}
\newtheorem{Def}[thm]{Definition}
\newtheorem{Cor}[thm]{Corollary}
\newtheorem{Ex}[thm]{Example}
\newtheorem{Pro}[thm]{Problem}
\newtheorem{Remark}[thm]{Remark}
\newtheorem{NN}[thm]{}
\renewcommand{\theequation}{e\,\arabic{equation}}

\newcommand{\for}{\,\,\, {\rm \,\,\,for\,\,\,}}
\newcommand{\rforal}{\,\,\,{\rm for\,\,\,all\,\,\,}}
\newcommand{\Ik}{ {\cal I}^{(k)}}
\newcommand{\Iz}{{\cal I}^{(0)}}
\newcommand{\Ii}{{\cal I}^{(1)}}
\newcommand{\Ip}{{\cal I}^{(2)}}

\begin{abstract}
Let $A$ be a unital separable amenable \CA\, and $C$ be a unital
\CA\, with certain infinite property. We show that two full
monomorphisms $h_1, h_2: A\to C$ are approximately unitarily
equivalent if and only if $[h_1]=[h_2]$ in $KL(A,C).$ Let $B$ be a
non-unital but $\sigma$-unital \CA\, for which $M(B)/B$ has
the certain infinite property. We prove that two full essential
extensions are approximately unitarily equivalent if and only if
they induce the same element in $KL(A, M(B)/B).$ The set of
approximately unitarily equivalence classes of full essential
extensions forms a group. If $A$ satisfies the Universal
Coefficient Theorem, it is can be identified with $KL(A, M(B)/B).$

\end{abstract}

\vspace{0.2in}


\section{Introduction}

The study of \CA\, extensions originated in the study of
essentially normal operators on the infinite dimensional separable
Hilbert space. The original Brown-Douglas-Fillmore theory gives a
classification of essential normal operators via certain Fredholm
related indices (see \cite{BDF1}). Later the
Brown-Douglas-Fillmore theory gives classification of essential
extensions of $C(X)$ by the compact operators (see \cite{BDF2},
\cite{Br2}). The study of \CA\,  extensions developed into
Kasparav's $KK$-theory and its application can be found not only
in operator theory and operator algebras  but also in differential
geometry and noncommutative geometry.

Let $0\to B\to E\to A\to 0$ be an essential extension of $A$ by
$B.$ The extension is determined by a monomorphism $\tau: A\to
M(B)/B,$ the Busby invariant. When $B$ is $\sigma$-unital stable
\CA\, then $KK^1(A,B)$ gives a complete classification of these
essential extensions --up to stable unitary equivalence. However,
$KK^1(A,B)$ gives  little information, if any, about unitary
equivalence classes of the above mentioned extensions when
$B\not={\cal K}$ in general. There are known examples in which
$KK^1(A,B)=\{0\}$ but inequivalent non-trivial extensions exist
(see Example 0.6 of \cite{Ln5}). There are also known examples in
which there are infinitely many inequivalent classes of trivial
extensions (see 7.4 and 7.5 of \cite{Ln4}). When $B$ is not
stable, $KK^1(A,B)$ certainly should not be used to understand
unitary equivalence classes of essential extensions mentioned
above.

There are a number of results in classification of essential
extensions ( up to unitary equivalence or approximate unitary
equivalence) when $B\not={\cal K}.$ Kirchberg's results
(\cite{K1}) on extensions in which $B$ is a non-unital purely
infinite simple \CA\, shows that $KK^1(A,B)$ can be used to
compute unitary equivalence classes of those extensions. When $B$
is a non-unital but $\sigma$-unital simple \CA\, with continuous
scale (see (6) below), then $M(B)/B$ is simple.  Classification of
essential extensions of a separable amenable \CA\, $A$ by $B$ (up to
approximate unitary equivalence) was obtained in \cite{Lnnext}
(for some special cases in which  $A=C(X),$ classification up to
unitary equivalence was obtained in \cite{Ln3}, \cite{Ln4} and
\cite{Lnamj}). In this case, $B$ may not be stable, therefore
$KK^1(A,B)$ is not used as invariant for essential extensions.
Results about extensions of AF-algebras  may be found in
(\cite{BE},\cite{GH} and \cite{EH}).

In this paper, we study full essential extensions. These are
essential  extensions $\tau: A\to M(B)/B$ so that $\tau(a)$ is a
full element for each nonzero element $a\in A.$ Since the Calkin
algebra $M({\cal K})/{\cal K}$ is simple, all essential extensions
by ${\cal K}$ are full. If $B$ is a non-unital but $\sigma$-unital
purely infinite simple \CA\, then $M(B)/B$ is also simple.
Therefore essential extensions by those \CA s are all full. The
homogeneuous extensions of $A$ by $C(X)\otimes {\cal K}$ studied
by Pimsner-Popa-Voiculescu (\cite{PPV1} and \cite{PPV2}) are all
full extensions. In all these three cases, $B$ is stable. There
are non-stable, non-unital but $\sigma$-unital \CA s which have
continuous scale. In that case essential extensions by these \CA s
are also full. Furthermore, if $A$ is a unital simple \CA\, and if
the monomorphism $\tau: A\to M(B)/B$ is unital, then the essential
extension induced by $\tau$ is always full  for any non-unital
\CA\, $B.$

With a  technical condition on $M(B)/B,$ we show that two full
essential extensions are approximately unitarily equivalent if
they induce the same element in $KL(A, M(B)/B)$ (see Theorem
\ref{MT1}) provided that $A$ is amenable and separable. When $A$
is assumed to satisfy the so-called (Approximate) Universal
Coefficient Theorem, we show that there is a bijective
correspondence between approximate unitary equivalence classes of
essential full extensions and $KL(A,M(B)/B).$ The advantage of
study these full extensions is that full extensions (in these
cases) are ``approximately absorbing". For stable $B,$ we show
that $KK^1(A,B)$ classifies the unitary equivalence classes of
essential full extensions. In this case, full extensions are
``purely large" in the sense of Elliott and Kucerovsky
(\cite{EK}).

The paper is organized as follows. Section 2 describes the main
results in this paper. Section 3 shows that for many stable \CA s
their corona algebras $M(B)/B$ satisfy the technical condition
(P1), (P2) and (P3). In Section 4, we show that there are examples
of non-stable, non-unital and $\sigma$-unital \CA s $B$ for which
$M(B)/B$ has the property (P1), (P2) and (P3). In Section 5, we
give  some modified results concerning  amenable \morp s. In
Section 6, we discuss certain commutants in the ultrapower of corona
algebras. In Section 7, we prove  Theorem \ref{MT1} mentioned
above. In  Section 8, we prove other main results described in
Section 2.

We will use the following convention:

(1) Ideals in this paper are always closed and two-sided.

(2) Let $A$ be a \CA\, and $p, q\in A$ be two projections.
We write $p\sim q$ if there exists $v\in A$ such that
$v^*v=p$ and $vv^*=q.$

(3) Let $A$ and $B$ be \CA s and $L_1, L_2: A\to B$ be linear
maps. Let ${\cal F}\subset A$ and $\ep>0,$ we write $L_1\sim_{\ep}
L_2$ on ${\cal F},$ if
$$
\|L_1(a)-L_2(a)\|<\ep\rforal a\in {\cal F}.
$$

(4) Let $A$ and $B$ be \CA s. A \morp\,
$L: A\to B$ is said to {\it amenable}, if for $\ep>0,$ there exists
an integer $n>0$ and two  \morp s $\phi: A\to M_n$ and
$\psi: M_n\to A$ such that
$$
\psi\circ \phi\sim_{\ep} L\,\,\,{\rm on}\,\,\,{\cal F}.
$$

(5) A \CA\, $A$ is said to be {\it  amenable} (or nuclear) if ${\rm id}_{A}$ is
amenable.

(6) Let $B$ be a non-unital but $\sigma$-unital simple
\CA. $B$ is said to have continuous scale, if there
exists an approximate identity $\{e_n\}$ of $B$ with
$e_{n+1}e_n=e_n$ such that, for each nonzero element $b\in B,$
there exists an integer $n>0$ for which
$e_{n+m}-e_n\lesssim b$ for all $m$ (see \cite{Lncs2}).

Let $e\in B$ be a nonzero projection and $T_e(B)$ be the set of
all traces $t$  on $B$ for which $t(e)=1.$ Let $B$ be a separable
non-unital simple \CA\, with real rank zero, stable rank one and
weakly unperforatated $K_0(B).$ If $\sup_n\{t(e_n)\}$ is a
continuous function on $T_e(B),$ then $B$ has continuous scale.

(7) Let $\{A_n\}$ be a sequence of \CA s. Denote by $c_0(\{A_n\})$
the ($C^*$-) direct sum of $\{A_n\}$ and denote by
$l^{\infty}(\{A_n\})$ the ($C^*$-) product of $\{A_n\}.$  We use
$q_{\infty}(\{A_n\})$ for the quotient
$l^{\infty}(\{A_n\})/c_0(\{A_n\}).$  When $A=A_n$ for all $n,$ we
write $c_0(A),$ $l^{\infty}(A)$ and $q_{\infty}(A)$ for
simplicity.

(8)
For each integer $n>0,$ define  $f_n\in C_0((0, \infty))$ as
follows
\begin{eqnarray}\label{fn}
f_n(t)=\begin{cases} 1 &\text{ if $t\ge 1/n$;}\\
                     {\rm linear} & \text{if $1/(n+1)\le t<1/n$;}\\
                     0 & \text{if $0\le t<1/(n+1)$}.
\end{cases}
\end{eqnarray}

(9) An element $a$ in a \CA\,$A$ is said be {\it full}, if
the ideal generated by $a$ is $A$ itself.
Let $A$ and $B$ be two \CA s and let $h: A\to B$ be a monomorphism.
The monomorphism $h$ is said to be {\it full}\,
if $h(a)$ is full for every nonzero $a\in A.$

(10) Let $a\in A_+$ be a nonzero element, we write $Her(a)$ for
the hereditary \SCA\, $\overline{aAa}$ generated by $a.$

{\bf Acknowledgement} This work started in summer 2003 when the
author was visiting East China Normal University. It is
partially supported  by National Science Foundation of U.S.A.

\section{Main results}

\begin{Def}\label{DP1}
{\rm Let $B$ be a unital \CA. We say that $B$ has property (P1) if
for every full element $b\in B$ there exist $x, y\in B$ such that
$xby=1.$  If $b$ is positive, it is easy to see that $xby=1$
implies that there is $z\in B$ such that $z^*bz=1.$}
\end{Def}

   It is obvious that an element $b$ is full if and only if $b^*b$
   is full.  It follows that $B$ has property (P1) if and only
   if for every full element $0\le b\le 1,$ there exists
   $x\in B$ such that $x^*bx=1.$

  Every unital purely infinite simple \CA\, has the property (P1).

It turns out that many other unital \CA s have the property (P1).
Let $A$ be a unital \CA\, and $B=A\otimes {\cal K}.$ In next
section we will show that $M(B)$ and $M(B)/B$ have property (P1)
for many such that $B$. In Section 3,  we will show that, for some
non-stable (but $\sigma$-unital)
 \CA\, $C,$  $M(C)$ and $M(C)/C$ may also have
property (P1).

\begin{Def}\label{Dp2}
{\rm Let $B$ be a unital \CA. We say that $B$ has property (P2),
if $1$ is proper infinite, i.e., there is a projection $p\not=1$ and
partial  isometries $w_1,w_2\in B$
such that $w_1^*w_1=1,$ $w_1w_1^*=p,$ $w_2^*w_2=1$ and $w_2w_2^*\le 1-p.$

It is easy to see that, for each
integer $n\ge 2$ and there are mutually
orthogonal and mutually equivalent projections
$s_{11},s_{22},...,s_{nn}$ such that $1_B\ge \sum_{i=1}^ns_{ii}$
and there exists an isometry $Z\in B$ such that $Z^*Z=1_B$ and
$ZZ^*=s_{11}.$ Let $C=s_{11}Bs_{11}.$ Then we may write
$M_n(C)\subset B.$
 }
\end{Def}

 It is clear that if $B$ is stable then $M(B)$ and $M(B)/B$ have
 property (P2).
\begin{Prop}\label{PM1}
Let $B$ be a unital \CA\, which has property (P1). Suppose that
$B$ contains two mutually orthogonal full elements. Then $B$ has
property (P2).
\end{Prop}

\begin{proof}
Let $0\le a, b\le 1$ be two mutually orthogonal full elements in
$B.$ Since $B$ has property (P1), there are $x,y\in B$ such that
$x^*ax=1$ and $y^*by=1.$ Let $v_1=a^{1/2}x$ and $v_2=b^{1/2}y.$
Then $v_i^*v_i=1$ and  $s_{11}=v_1v_1^*$ and  $s_{22}=v_2v_2^*$ are two
projections. Thus $B$ has property (P2).
\end{proof}

Every purely infinite \CA\, (not necessary simple; see \cite{KR})
has property (P1) and (P2).

\begin{Def}\label{Dp3}
{\rm Let $B$ be a unital \CA. We say that $B$ has property (P3),
if for any separable  \SCA\, $A\subset B,$ there exists a sequence
of sequences of  elements $\{\{a_n^{(i)}\}: i=1,2,...\}$ in $B$
with $0\le a_n\le 1$ such that
$$
\lim_{n\to\infty}\|a_n^{(i)}c-ca_n^{(i)}\|=0\rforal c\in
A,\,i=1,2,...,
$$
$\lim_{n\to\infty}\|a_n^{(i)}a_n^{(j)}\|=0$ if $i\not=j$ and for
each $i,$ and  $\{a_n^{(i)}\}$ is a full element in $l^{\infty}(B ).$

}
\end{Def}

Even though property (P3) looks more complicated than (P1) and (P2),
it will be shown (see \ref{IIPP3} below) that $M(B)/B$ has property
(P3) for all $B=C\otimes {\cal K},$ where $C$ is a unital \CA\,
and for all $B$ which have continuous scale and for many other
non-unital $\sigma$-unital \CA s $B.$

\vspace{0.1in}

\begin{Prop}\label{MPoinf}
Let $B=C\otimes C_1,$ where $C_1$ is a unital separable  amenable
purely infinite simple \CA.
Then $B$ has property (P1), (P2) and (P3).
\end{Prop}

Let $B$ be a non-unital but $\sigma$-unital \CA\, and $A$  be a
unital separable amenable \CA. We study essential extensions of
the following form:
\begin{eqnarray}\label{IVext}
0\to B\to E\to A\to 0.
\end{eqnarray}
Using the Busby invariant, we study monomorphisms $\tau: A\to
M(B)/B.$ We will only consider the case that the corona algebra
$M(B)/B$ has the property (P1),(P2) and (P3).

\begin{Def}
{\rm An essential extension $\tau: A\to M(B)/B$ is said to be {\it full},
if $\tau$ is a full monomorphism
. An extension $\tau$ is {\it weakly
unital} if $\tau$ is unital monomorphism.
If $A$ is a unital simple \CA\, then every weakly unital essential
extension is full. If $M(B)/B$ is simple, then every essential extension
is full.}
\end{Def}

\begin{Def}
{\rm Let $A$ be a unital separable \CA\, and $C$ be a unital \CA.
Suppose that $h_1, h_2: A\to C$ are two \hm s. We say $h_1$ and
$h_2$ are {\it approximately unitarily equivalent} if there exists
a sequence of partial isometries  $u_n\in C$ such that
$u_n^*h_1(1_A)u_n=h_2(1_A),$ $u_nh_2(1_A)u_n^*=h_1(1_A)$ and
$$
\lim_{n\to\infty}\|{\rm ad}\, u_n\circ h_1(a)-h_2(a)\|=0\rforal
a\in A.
$$
Note that if both $h_1$ and $h_2$ are
unital, $u_n$ can be chosedn to be unitaries.

Let $B$ be a non-unital but $\sigma$-unital \CA. Two essential
extensions of $A$ by $B$ are said to be approximately unitarily
equivalent if the corresponding Busby invariants $\tau_1, \tau_2:
A\to M(B)/B$ are approximately unitarily equivalent.

Recall that $\tau: A\to M(B)/B$ is trivial if there is a
monomorphism $h: A\to M(B)$ such that $\pi\circ h=\tau,$ where
$\pi:M(B)\to M(B)/B$ is a quotient map.
 In the case that
$B=C\otimes {\cal K},$ where $C$ is a $\sigma$-unital \CA\,
$\tau_1$ and $\tau_2$ are {\it stably} unitarily equivalent if
there exists a trivial extension $\tau_0: A\to M(B)/B$ and a
unitary $u\in M_2(M(B)/B)$ such that ${\rm ad} u\circ
(\tau_1\oplus \tau_0)=\tau_2\oplus \tau_0.$

 }
\end{Def}
 Let ${\bf Ext}(A,B)$ be the {\it stable} unitary equivalence
 classes of extensions of the form (\ref{IVext}). When
 $A$ is a separable amenable \CA\, ${\bf Ext}(A,B)$ may be
 identified with $KK^1(A,B).$ When $A$ satisfies the Universal
 Coefficient Theorem, $KK^1(A,B)$ may be computable.
 However, as mentioned in the introduction, $KK^1(A,B)$ may not
 provide any useful information about  unitary equivalence
 of extensions in general. In particular, when $B$ is not stable,
 $KK^1(A,B)$ should not be used to describe unitary equivalence
 classes
 of essential extensions.

 The first main result of this paper is the following:

 \begin{thm}\label{MT1}
Let $A$ be a unital separable amenable  \CA\, and  $B$ be a
non-unital but $\sigma$-unital \CA\, so that $M(B)/B$  has the
property (P1), (P2) and (P3).  Suppose that $\tau_1, \tau_2: A\to
M(B)/B$ are two full monomorphisms. Then $\tau_1$ and $\tau_2$ are
approximately unitarily equivalent if and only if
$$
[\tau_1]=[\tau_2]\,\,\,\,{\rm in}\,\,\, KL(A, M(B)/B).
$$
\end{thm}

We will describe $KL(A, C)$ in \ref{DKL}.
Theorem \ref{MT1} is an easy corollary of the following
theorem.

\begin{thm}\label{MTad}
Let $A$ be a unital separable amenable  \CA\, and  $B$ be a unital
\CA\, which has property (P1), (P2) and (P3). Suppose that $h_1,
h_2: A\to B$ are  two full monomorphisms. Then $h_1$ and $h_2$ are
approximately unitarily equivalent, i.e, there exists a sequence
of partial isometries  $u_n\in B$ such that $u_n^*u_n=h_1(1_A),$
$u_nu_n^*=h_2(1_A)$ and
$$
\lim_{n\to\infty} {\rm ad}\, u_n\circ h_1(a)=h_2(a)\rforal a\in A
$$
if and only if $[h_1]=[h_2]$ in $KL(A,B).$
\end{thm}

\begin{Cor}\label{MC1}
Let $A$ be a unital separable amenable simple \CA\, and $B$ be a
non-unital but $\sigma$-unital \CA\, so that $M(B)/B$  has the
property (P1), (P2) and (P3).  Suppose that $\tau_1, \tau_2: A\to
M(B)/B$ are two weakly unital essential extensions. Then $\tau_1$
and $\tau_2$ are approximately unitarily equivalent if and only if
$$
[\tau_1]=[\tau_2]\,\,\,\,{\rm in}\,\,\, KL(A, M(B)/B).
$$
\end{Cor}

\begin{Def}\label{MDg}
Let $A$ be a unital separable amenable \CA\, and $B$ be a unital
\CA\, which has property (P2). Fix a full monomorphism $j_o: A\to
{\cal O}_2\to B.$ Note (P2) implies such full monomorphisms do
exist.  Let $h_1, h_2: A\to B\otimes {\cal K}$ be two \hm s. We
write $h_1\sim h_2$ if $h_1\oplus j_o$ is approximately unitarily
equivalent to $h_2\oplus j_o.$ Denote by $H(A,B)$ be ``$\sim$ "
equivalent classes of those \hm s.
\end{Def}

\begin{Prop}\label{MPg1}
Let $A$ be a unital separable amenable \CA\, and $B$ be a unital
\CA\, which has property (P2). Then $H(A,B)$ is a group with the
zero element $[j_o].$
\end{Prop}

\begin{Cor}\label{MPg2}
Let $A$ be a unital separable amenable \CA\, and $B$ be a unital
\CA\, which has property (P1), (P2) and (P3). Let $H_{f}(A,B)$ be
the approximate unitary equivalence classes of full monomorphisms
from $A$ to $B\otimes {\cal K}.$ Then $H_f(A,B)$ is a group with
the zero element $[j_o].$
\end{Cor}

\begin{Def}\label{MDa1}
{\rm Let $A$ be a unital separable amenable \CA\, and $B$ be a
non-unital but $\sigma$-unital \CA. Denote by ${\bf
Ext}^f_{ap}(A,B)$ the approximate unitary equivalence classes of
full essential extensions. Denote by $\tau_o: A\to M(B)/B$ an
essential extension which factors through ${\cal O}_2.$ Note that
$[\tau_o]=0$ in $KL(A, M(B)/B).$ Suppose that $M(B)/B$ has
property (P1), (P2) and (P3). It follows this $\tau_o$ is unique
up to approximately unitary equivalence,
 by
\ref{MC1}.
Let $\tau_1, \tau_2: A\to
M(B)/B$ be two essential full extensions. Since $M(B)/B$ has
property (P2), there are partial isometries $z_1, z_2\in M(B)/B$
such that $z_1^*z_1=1_{M(B)/B},$ $z_1z_1^*=\tau_1(1_A),$
$z_2^*z_2=1_{M(B)/B}$ and $z_2z_2^*=\tau_2(1_A).$ Define
$[\tau_1]+[\tau_2]=[{\rm ad}\,z_1\circ \tau_1\oplus {\rm
ad}\,z_2\circ \tau_2].$

Note this is well defined, since $[\tau_o]=0$ in $KL(A. M(B)/B)$
and ${\rm ad}\, z_1\circ \tau \oplus{\rm ad}\,z_2\circ \tau_o$ is
approximately unitarily equivalent to $\tau$ by \ref{MC1}.
 With this addition ${\bf Ext}_{ap}^f(A,B)$
forms a semigroup.
By \ref{MPg2}, we have the following.}

\end{Def}

\begin{Cor}\label{MPg3}
Let $A$ be a unital separable amenable \CA\, and $B$ be a
non-unital but $\sigma$-unital \CA\, for which $M(B)/B$ has
property (P1), (P2) and (P3). Then ${\bf Ext}^f_{ap}(A,B)$ is a
group with zero element $[\tau_o],$ where $\tau_o: A\to M(B)/B$ is
a full monomorphism which factors through ${\cal O}_2.$
\end{Cor}

If furthermore, $A$ satisfies so-called Approximate Universal
Coefficient Theorem (AUCT) (see \ref{DKL} below), then we have the
following.

\begin{thm}\label{MT2}
Let $A$ be a unital separable amenable \CA\, which satisfies
Approximate Universal Coefficient Theorem and $B$ be a non-unital
but $\sigma$-unital \CA\, so that $M(B)/B$  has the property (P1),
(P2) and (P3).  Then there is a bijection $\Gamma$ from ${\bf
Ext}_{ap}^f(A,B)$ onto $KL(A, M(B)/B).$
\end{thm}

Approximate Universal Coefficient Theorem will be briefly
discussed in \ref{DKL} and \ref{DKunder}.
It should be noted that, when $B$ is not stable, $K_i(M(B)/B)$ is
very different from $K_i(SB),$ $i=0,1.$ (see 1.7 of \cite{Lnnext}).

 In the
special case that $B=C\otimes {\cal K},$ where $C$ is a unital
\CA, we have the following theorem:

\begin{thm}\label{MT3}
Let $A$ be a unital separable amenable  \CA\, and $B=C\otimes
{\cal K},$ where $C$ is a unital \CA\, so that  $M(B)/B$ has the
property (P1).  Suppose that $\tau_1, \tau_2: A\to M(B)/B$ are two
full essential extensions. Then $\tau_1$ and $\tau_2$ are
unitarily equivalent if and only if
$$
[\tau_1]=[\tau_2]\,\,\,\,{\rm in}\,\,\, KK^1(A,B).
$$

Moreover, if $x\in KK^1(A,B),$ then there is a full essential
extension $\tau: A\to M(B)/B$ such that $[\tau]=x.$
\end{thm}

\begin{thm}\label{MT4}
Let $A$ be a unital separable
amenable  \CA\, and $B=C\otimes {\cal K},$
where $C$ is a unital \CA\, for which
the tracial state space $T(C)\not=\emptyset.$
Suppose that there is $d>0$ for which $C$ satisfies the following:

{\rm (1)}
   if $p, q\in B$ are two projections then  $t(p)>d+t(q)$ for all $t\in T(C) $
   implies $q\sim p$ in $B;$

 {\rm   (2)} if $1\ge b\ge 0$ in $M_k(C)$ such that $\tau(b)>\alpha+d$ for
   all $\tau\in T(A),$ then there is a projection $e\in
   \overline{bM_k(A)b}$ such that $\tau(e)>\alpha$ for all $\tau\in
   T(A).$

Then two essential full extensions $\tau_1, \tau_2: A\to M(B)/B$
are unitarily equivalent if and only if
$$
[\tau_1]=[\tau_2].
$$
\end{thm}

\begin{Remark}
{\rm  In the case that $B={\cal K},$ Theorem \ref{MT3} is  the
classical Brown-Douglas-Fillmore  theorem. Note in this case,
$M({\cal K})/{\cal K}$ is a purely infinite simple \CA. It has
property (P1) (as well as (P2) and (P3)) and every essential
extension is full. Let $X$ be a compact metric space with finite
dimension $d.$ When $B=C(X)\otimes {\cal K},$ $M(B)/B$ has
property (P1) (see \ref{CI4}). Theorem \ref{MT3} (or \ref{MT4})
deals with the extensions studied by Pimsner-Popa-Voiculescu (see
\cite{PPV1} and \cite{PPV2}). When $B$ is a non-unital purely
infinite simple \CA\, this is obtained by Kirchberg. This theorem
is closely  related to a result of Elloitt and Kucerovsky
(\cite{EK}), see \ref{Rf} for a discussion.}

\end{Remark}

\section{\CA s have property (P1), (P2) and (P3)}

 Let
$A$ be a unital \CA. Denote by $T(A)$ (or $T$ if no confusion
exits) the set of tracial states on $A.$ If $t\in T(A),$ we extend
$t$ to a trace ($t\otimes Tr$) on $A\otimes M_n$ by defining
$t((a_{ij})=\sum_{i=1}^n t(a_{ii}).$ We further use $t$ for the
trace defined on a dense set on $A\otimes {\cal K}.$ If $a\in
A\otimes {\cal K}_+,$ then $t(a)$ is well defined (although it could be
infinity). Suppose that $h_n\in A\otimes {\cal K}_+$ such
that $h_n\nearrow h\in A^{**}.$ Then one has
$t(h)=\lim_{n\to\infty} t(h_n).$ These conventions will be used in
this section.

The following lemma is certainly known

\begin{Lem}\label{Ifull}
Let $A$ be a unital \CA\, and $I$ be a $\sigma$-unital ideal of
$A.$ If $a\in (A/I)_+$ is a full element, then there exists a full
element $b\in A_+$ such that $\pi(b)=a,$ where $\pi: A\to A/I$ is
the quotient map.
\end{Lem}

\begin{proof}
Since $a\in (A/I)_+$ is full,  there are $x_1,x_2,...,x_m\in A/I$
such that
\begin{eqnarray}
\sum_{i=1}^mx_i^*ax_i=1.
\end{eqnarray}
It follows that there are $c\in A_+$ and $y_1,y_2,...,y_m\in A$
such that $\pi(c)=a$ and $1-\sum_{i=1}^m y_i^*cy_i\in I.$ Let $e$ be a strictly
positive element of $I.$ Put $b=c+e.$ Denote by $J$ the ideal
generated by $b.$ Since $b\ge c$ and $b\ge e,$ both $c$ and $e$
are in $J.$ It follows that $I\subset J.$ Since $\sum_{i=1}^m
y_i^*cy_i\in J,$ it follows that $1\in J.$ Thus $J=A.$ Therefore
$b$ is full.
\end{proof}

\begin{Cor}\label{ICid}
Let $A$ be a unital \CA\, and $I$ be a $\sigma$-unital ideal of
$A.$ If $A$ has property (P1), then so does $A/I.$
\end{Cor}

\begin{Lem}\label{Ifull2}
Let $A$ be a unital \CA\, $B=A\otimes {\cal K},$ Suppose that
$a\in M(B)$ is an element for which $b=\pi(a)$ is full in
$M(B)/B,$ where $\pi: M(B)\to M(B)/B$ is the quotient map. Then
$a$ is full in $M(B).$
Furthermore, $M(B)/B$ has property (P1) so does $M(B).$
\end{Lem}

\begin{proof}
There are $x_1,x_2,...,x_m\, y_1,...,y_m \in M(B)/B$ such that
$\sum_{i=1}^m x_i b y_i=1.$ Then there are $w_1,w_2,...,w_m,$
$z_1,z_2,...,z_m\in M(B)$ such that
$$
1-\sum_{i=1}^m w_iaz_i\in B.
$$
Let $\{e_{ij}\}$ be a system of matrix units. Put
$E_n=\sum_{i=1}^ne_{ii}.$ Then $\{E_n\}$ is an approximate
identity consisting of projections. It follows that there exits
$n>0$ such that
\begin{eqnarray}
\|\sum_{i=1}^m(1-E_n)w_iaz_i(1-E_n)-(1-E_n)\|<1/2
\end{eqnarray}
Thus there exists $s\in (1-E_n)M(B)(1-E_n)$ such that
\begin{eqnarray}
\sum_{i=1}^ms^*(1-E_n)w_iaz_i(1-E_n)s=1-E_n.
\end{eqnarray}
But there exists $V\in M(B)$ such that $V^*(1-E_n)V=1.$ Therefore
$a$ is full.

For the last statement, we take $m=1$ in the above argument.
\end{proof}

\begin{Prop}
Let $B$ be a unital purely infinite simple \CA. Then $M(B\otimes
{\cal K})$ and $M(B\otimes {\cal K})/A\otimes {\cal K}$ have the
property (P1).
\end{Prop}

\begin{proof}
It follows from \cite{Zh1} that $M(B\otimes {\cal K})/B\otimes
{\cal K}$ is purely infinite and simple. Therefore $M(B\otimes
{\cal K})/B\otimes {\cal K}$ has the property (P1).
It follows from \ref{Ifull2} that $M(B\otimes {\cal K})$ has
property (P1).
\end{proof}

\begin{thm}\label{TI1}
Let $B=A\otimes {\cal K},$ where $A$ is a unital separable  \CA\,
for which $T(A)\not=\emptyset.$ Let $d>0.$
Suppose $A$  satisfies the following:

{\rm (1) } if $p, q\in B$ are two projections then  $t(p)>d+t(q)$
for all $t\in T(A) $ implies $q\sim p$ in $B;$

{\rm (2)} if $1\ge b\ge 0$ in $M_k(A)$ such that
$\tau(b)>\alpha+d$ for all $\tau\in T(A)$ (and some $\alpha>0$),
then there is a projection $e\in \overline{bM_k(A)b}$ such that
$\tau(e)>\alpha$ for all $\tau\in T(A).$

Then $M(B)$ and $M(B)/B$ have  property {\rm (}P1{\rm )}.
\end{thm}

\begin{proof}
Let $b\in M(B)$ be a full element.
Without loss of generality,
we may assume that $0\le b\le 1.$
Let $\{e_{ij}\}$ be the system of matrix unit for ${\cal K}$ and
$E_n=\sum_{k=1}^ne_{ii}.$
It follows that $E_nbE_n$ converges to $b$ in the strict
topology. Furthermore $b^{1/2}E_nb^{1/2}$ increasingly
converges to $b$ in the strict topology.

Since $b$ is full, there are $x_1,x_2,...,x_m\in M(B)$ such that
$$
\sum_{k=1}^mx_i^*bx_i=1.
$$
Let $\tau\in T(A)$ be a tracial state. We extend $\tau$ to $B_+$
and then $M(B)_+$ in a usual way. Let $T$ be the set of all
(densely defined ) traces on $M(B)_+$ whose
restrictions  to
$A$ are  tracial states. With the usual weak *-topology, $T$ is a
compact convex set.

Because $b^{1/2}x_i^*x_ib^{1/2}\le \|x_i\|^2b,$ one has
$$
\tau(x_i^*bx_i)=\tau(b^{1/2}x_i^*x_ib^{1/2})\le \|x_i\|^2\tau(b)
$$
for all $\tau\in T(A)$ Therefore
$$
\sum_{i=1}^m \tau(x_i^*bx_i)\le (\sum_{i=1}^m\|x_i\|^2)\tau(b)
$$
for all $\tau\in T(A).$ Since $\tau(1)=\infty,$ it follows that
$\tau(b)=\infty.$ Because $b^{1/2}E_nb^{1/2}\nearrow b,$ and
because $T$ is compact, by the Dini's theorem,
$\tau(b^{1/2}E_nb^{1/2})\to\infty$ uniformly on $T.$ Since
$\tau(E_nbE_n)=\tau(b^{1/2}E_nb^{1/2})$ for all $\tau\in T,$
$\tau(E_nbE_n)\nearrow \infty$ uniformly on $T.$ There is $n(1)\ge
1$ such that
$$
\tau(E_{n(1)}bE_{n(1)})>1 +2d \,\,\,{\rm
for\,\,\,all}\,\,\,\tau\in T.
$$
Let $A_1$ be the hereditary \SCA\, of $B$ generated by $E_{n(1)}bE_{n(1)}.$
It follows from assumption (2) that there is a projection $p_1\in
A_1$ such that $\tau(p_1)>1+d$ for all $\tau\in T.$ It follows
that there is $v_1\in B$ such that $v_1^*v_1\le p_1$ and
$v_1v_1^*=E_1.$ There are  non-negative continuous function
$f,g\in C_0((0,2\|b\|]$ such that $gf=f$ and
$$
\|f(E_{n(1)}bE_{n(1)})p_1f(E_{n(1)}bE_{n(1)})-p_1\|<1/4.
$$
It  follows (see A8 \cite{Eff}) that there is a projection
$q_1\in {\overline{f(E_{n(1)}bE_{n(1)})Bf(E_{n(1)}bE_{n(1)})}}$
such that $q_1$ is unitarily equivalent to $p_1.$
Since $gf=f,$ we conclude that $gq_1=q_1.$
By functional calculus, we see that there are $f_1\in A_1$ such that
$$
f_1E_{n(1)}bE_{n(1)}f_1=g.
$$
Thus we obtain $z_1\in E_{n(1)}BE_{n(1)}$ such that
$$
z_1^*bz_1=z_1^*E_{n(1)}bE_{n(1)}z_1=E_1.
$$
Note that $\tau((1-E_{n(1)})bE_{n(1)})=\tau(bE_{n(1)}(1-E_{n(1)}))=0.$
It follows that
$$
\tau((1-E_{n(1)})b(1-E_{n(1)}))=\tau((1-E_{n(1)})b).
$$
Since $\tau(E_{n(1)}bE_{n(1)})<\infty,$
for all $\tau\in T,$ we conclude that
$$
\tau((1-E_{n(1)})b(1-E_{n(1)})=\infty\,\,\,\,\,{\rm for\,\,\,all}\,\,\,
\tau\in T.
$$
From the above argument, we obtain $n(2)>n(1)$ and
$z_2\in (E_{n(2)}-E_{n(1)})B(E_{n(2)}-E_{n(1)})$ such that
$$
z_2^*bz_2=z_2^*(E_{n(2)}-E_{n(1)})b(E_{n(2)}-E_{n(1)})z_2=E_2-E_1.
$$
Continuing this fashion, we obtain a sequence $\{n(k)\}$ with
$n(k+1)>n(k)$ and\\
 $z_k\in
(E_{n(k+1)}-E_{n(k)})B(E_{n(k+1)}-E_{n(k)})$ such that
$$
z_k^*bz_k^*=z_k^*(E_{n(k+1)}-E_{n(k)})b(E_{n(k+1)}-E_{n(k)})z_k=
E_{k+1}-E_k,
$$
$k=1,2,....$ It follows that $z=\sum_{k=1}^{\infty} z_k\in M(B)$
since the sum converges in the strict topology. Furthermore
we have
$$
z^*bz=1.
$$
This shows that $M(B)$ has the property (P1).

By \ref{ICid}, $M(B)/B$ also has property (P1).
\end{proof}

From \ref{TI1}, we have the following corollaries.

\begin{Cor}\label{CI1} Let $A$ be a unital AF-algebra and
$B=A\otimes {\cal K}.$ Then $M(B)$ and $M(B)/B$ have  property
(P1).
\end{Cor}

\begin{proof}
Clearly $A$ satisfies (1) in \ref{TI1} with any $d>0.$
To see  that $A$ satisfies (2), we let $1\ge b\ge 0$
be an element in $M_n(A)$ such that $\tau(b)>\alpha+d$
for all $\tau\in T.$ Let $C=\overline{bM_n(A)b}$ and let
$\{e_n\}$ be an approximate identity for $C$ consisting of projections.
Then
$\|e_nbe_n-b\|\to 0$ as $n\to\infty.$ Since $0\le b\le 1,$ it follows
that $\tau(e_n)>\alpha+d$ for some $n>0$ and all $\tau\in T.$
\end{proof}

The proof of the corollary implies the following:

\begin{Cor}\label{CI2}
Let $A$ be a unital separable \CA\, for which
$T(A)\not=\emptyset$  and which satisfies (1) in \ref{TI1} and has real rank
zero. Then $M(B)$ and $M(B)/B$ have property (P1), where
$B=A\otimes {\cal K}.$
\end{Cor}

 \begin{Cor}\label{CI3}
   Let $B=A\otimes {\cal K},$ where $A$ is a unital simple  \CA\,
   with real rank zero,
   stable rank one and weakly unperforated $K_0(A)$
   Then both $M(B)$ and $M(B)/B$ have the property (P1).
   \end{Cor}

\begin{Cor}\label{CI4}
Let $A=C(X),$ where $X$ is a compact Hausdorff space with
covering dimension $d.$
Then $M(A\otimes {\cal K})$ and $M(A\otimes {\cal K})/A\otimes {\cal K}$
have property (P1).
\end{Cor}

\begin{proof}
Suppose that $e, f\in A\otimes {\cal K}$ are two projections. It
is clear that we may assume that $e, f\in M_n(C(X))$ for some
integer $n>0.$ Suppose that $\tau(e)>\tau(f)+d+1$ for all $t\in
T(A).$ It follows that for each $x\in X,$ the rank of $e(x)$ is
greater than $d+1$ + the rank of $f(x).$ It follows from 8.1.2 and
8.1.6 in \cite{Hu} (see 6.10.3 (d) of \cite{B}) that $f\lesssim
e.$ So (1) in \ref{TI1} holds (for $(d+1)/2$).

For (2), let $1\ge b\ge 0$ be an element in $M_k(C(X))$ for which
$\tau(b)>\alpha+(d+1).$
Let  $f_n$ be as in (\ref{fn}).
It follows that for some large $n,$ $\tau(f_n(b))>\alpha+(d+1)$
for all $\tau\in T(A).$ Thus, for each $\xi\in X,$ the rank of
$f_n(b)(\xi)$ is at least $\alpha+(d+1).$ By Lemma C in
\cite{BDR}, there is a projection $e\in \overline{bM_k(A)b}$ such
that the rank of $e(\xi)$ is greater $\alpha$ for all  $\xi\in X.$
It follows that $\tau(e)>\alpha$ for all $\tau\in T(A).$
\end{proof}

To discuss property (P2), we  begin with the following
easy observation.

\begin{Prop}\label{IIP2}
Let $B$ be a unital \CA\, which has property (P2). Then, for any
integer $n>0,$ there are $s_{11},s_{22},...,s_{nn}$ such that
$1_B\ge \sum_{i=1}^n s_{ii}$ and there exists an isometry $Z\in B$
such that $ZZ^*=e_{11}.$ Moreover

{\rm (1)} if for some $n\ge 2,$ $1_B=\sum_{i=1}^ns_{ii},$ then
there exists a unital embedding from ${\cal O}_n$ to $B;$

{\rm (2) } there is a unital embedding from ${\cal O}_{\infty}$ to
$B$ and

{\rm (3)} there  exists a full embedding $j: {\cal O}_2\to
B.$

Conversely, if there is a unital  embedding from ${\cal
O}_{\infty}$ to $B,$ then $B$ has property (P2).
Furthermore, if $B$ admits a full embedding from ${\cal O}_2,$
then $B$ has property (P2).
\end{Prop}

\begin{Prop}

{\rm (1) } Let $A$ be a unital \CA\, and $B=A\otimes {\cal K}.$
Then $M(B)$ and $M(B)/B$ has property (P2).

{\rm (2)} Let $A$ be a non-unital $\sigma$-unital simple \CA\,
which has continuous scale. Then $M(A)/A$ has property (P2)

{\rm (3)} Let $A$ be a unital purely infinite simple \CA\, and
$B=C_0(X, A),$ where $X$ is a locally compact Hausdorff space.
Then $M(B)$ and $M(B)/B$ have property (P2).

\end{Prop}
\begin{proof}
For (3), we note there is a unital embedding from ${\cal
O}_{\infty}$ to $A$ and the constant maps from $X$ into $A$ are
in $C^b(X,A)=M(B).$
\end{proof}

Now we will turn to property (P3). Every unital purely infinite
simple \CA\, has property (P3). This follows from 2.6 of
\cite{Lnnext}. Therefore, if $B$ is a
non-unital but $\sigma$-unital simple \CA\, with continuous scale,
then $M(B)/B$ has property (P3).

\begin{Prop}\label{PI2}
Let $B$ be a unital \CA\, which has the property (P1). Suppose
that $0\le a, b\le 1$ where $ab=a$ and $a$ is full. Then there
exists $x\in B$ with $\|x\|\le 1$ such that
\begin{equation}
x^*bx=1.
\end{equation}
\end{Prop}

Note that the proposition includes the case that $a$ is a full
projection.

\begin{proof}
There is $z\in B$ such that $z^*az=1.$ Then $a^{1/2}zz^*a^{1/2}=p$
must be a projection. Moreover, $p\in Her(a).$ Therefore $pb=p.$
Put $v=a^{1/2}z.$ Then $v^*v=1$ and $vv^*=p.$ In particular,
$\|v\|=1.$ Now
\begin{eqnarray}
  1\ge \|b\|v^*v\ge  v^*bv\ge v^*pv=1.
\end{eqnarray}
We conclude that $v^*bv=1.$
\end{proof}

\begin{Prop}\label{IIPP3}
Let $A$ be a unital \CA, and $B=A\otimes {\cal K}.$
Then $M(B)/B$ has property (P3).
\end{Prop}
\begin{proof}
Let $\pi: M(B)\to M(B)/B$ be the quotient map  and $D$ be a
separable \CA. Let
$\{e_{i,j}\}$ be a system of matrix unit for ${\cal K}.$ Denote by
$E_n=\sum_{i=1}^ne_{i,i}.$ It is known (see
3.12.14 of \cite{Ped} and the proof of 5.5.3 of \cite{Lnb}) that there are
$\{e_n\}\subset {\rm Conv}\{E_n: n=1,2,...,\}$ such that
\begin{eqnarray}\label{IILcom1}
e_{n+1}e_n=e_n\andeqn \|e_na-ae_n\|\to 0,\,\,\, {\rm as}\,\,\, n\to\infty
\end{eqnarray}
for all $a\in D.$

Suppose that $e_n=\sum_{i=1}^{k(n)} \alpha_iE_i,$ where $\alpha_i$ are
non-negative scalars with $\sum_{i=1}^{k(n)}\alpha_i=1.$ There are
$0\le \beta_j\le 1$ such that $e_n=\sum_{i=1}^{k(n)} \beta_je_{jj}.$
Since, for each $i,$
\begin{eqnarray}
\|e_me_{ii}-e_{ii}\|\to 0\,\,\,\,\,\,{\rm as}\,\,\, m\to \infty
\end{eqnarray}
there is $N(n)>0$ such that, for each $m>N(n),$
$e_m=\sum_{i=1}^{k(m)} \beta_ie_{ii}$ with $\beta_{k(n)+1}>1/2.$ It
follows that $(e_m-e_n)e_{k(n)+1,k(n)+1}=\beta_{k(n)+1}e_{k(n)+1,k(n)+1}.$ By
passing to a subsequence if necessary, without loss of generality,
we may assume that $(e_{n+1}-e_n)e_{k(n)+1,k(n)+1} =\lambda_n
e_{k(n)+1,k(n)+1}$ for
some $\lambda_n>1/2.$ Now let $F\subset {\mathbb N}$ be an
infinite subset. Then
\begin{eqnarray}
b_F=\sum_{n\in F}(e_{n+1}-e_n)\ge (1/2)\sum_{n\in F} e_{k(n)+1,k(n)+1}.
\end{eqnarray}
It follows that $b_F$  is a full positive element in $M(B).$
Suppose that $\{F_n\}$ is a sequence of infinite subsets of
${\mathbb N}.$  Then, by \ref{PI2}, $\pi(\{\sum_{j\in
F_n}e_{k(j)+1,k(j)+1})\}$ is full in $l^{\infty}(M(B)/B).$  So
$\{\pi(b_{F_n})\}$ is full in $l^{\infty}(M(B)/B).$

By (\ref{IILcom1}), $\pi(b_F)$ commutes with $\pi(d)$ for each
$d\in D.$ Also  by (\ref{IILcom1}), if $|n-m|\ge 2,$
\begin{eqnarray}
(e_{n+1}-e_n)(e_{m+1}-e_m)=0.
\end{eqnarray}
It follows that $b_Fb_{F'}=0,$ if $|n-m|\ge 2$ for any $n\in F$
and any $m\in F'.$ Note that one may write $b_F=\sum_{n\in S(F)}
\lambda_n e_{n,n},$ where each $0<\lambda_n\le 1$ is a positive
number and $S(F)$ is an infinite subset of ${\mathbb N}.$

It is easy to find a family of (disjoint) infinite subsets
$\{F_{i,j}: i,j=1,2,...\}$ of ${\mathbb N}$ such
that $|n-m|\ge 2$ for any $n\in S_{i,j}$ and any $m\in S_{i',j'},$
if $i\not=i',$ or $j\not=j',$
Define $S_{i,j}=S(F_{i,j})$ as above. We note that $S_{i,j}\cap S_{i',j'}
=\emptyset,$ if $i\not=i'$ or $j\not=j'.$
Write $b_{i,j}$ for $b_{F_{i,j}}.$ It follows that $M(B)/B$ has
property (P3).
\end{proof}

\section{Non-stable cases}
In \cite{KR}, Kirchberg and R\o rdam extended the notion of purely
infinite \CA s to non-simple \CA s. Let $C_1$ be a unital \CA\,
and $C_2$ be a unital separable purely infinite simple \CA. Then
$C_1\otimes C_2$ is purely infinite (4.5 in \cite{KR}). Therefore,
for any unital \CA\, $C,$ $B=C\otimes {\cal O}_{\infty},$ has
property (P1) and (P2) as well as (P3).

\vspace{0.1in}
 {\bf Proof of Proposition \ref{MPoinf} }
\begin{proof}
By 4.5 in \cite{KR}, $B$ is purely infinite. It follows $B$ has
(P1) and (P2). Let $A$ be separable \SCA\, of $B.$ There is a
separable \SCA\, $C_0\subset C$ such that $A\subset C_0\otimes
C_1.$ It follows from \cite{KP} that $C_1\otimes {\cal O}_{\infty}
\cong C_1$ and  it follows from 7.2.6 of \cite{Rb} and 3.12 of
\cite{KP} that there is a sequence of unital monomorphisms
$\phi_n: {\cal O}_{\infty}\to C_0\otimes C_1$ such that
\begin{eqnarray}
\lim_{n\to\infty}\|\phi_n(x)a-a\phi_n(x)\|=0\rforal a\in C_0\otimes {\cal C}_1.
\end{eqnarray}
Let $\{e_k\}$ be a sequence of nonzero mutually orthogonal
projections in ${\cal O}_{\infty}.$ Define
$a_n^{(i)}=\phi_n(e_i),$ $n,i=1,2,....$ One checks that
$a_n^{(i)}$ satisfies the requirements in \ref{Dp3}.
\end{proof}

There are $\sigma$-unital but {\it non-stable} separable \CA s $B$
for which the corona \CA\,  $M(B)/B$ has  property (P1), (P2) as
well as (P3). For example, when $B$ has continuous scale (see
\cite{Lncs1} and \cite{Lncs2}) $M(B)/B$ is a purely infinite
simple \CA\, (see \cite{Lncs2}). So in those cases $M(B)/B$ has
property (P1), (P2) as well as (P3). There are other non-stable
separable \CA s $B$ for which  $B$ has property (P1), (P2) and
(P3).

To make a point,  we will present
a very simple  example of
non-stable  $\sigma$-unital \CA\, $B$ for which $M(B)/B$ is not
simple but both $M(B)$ and $M(B)/B$ have property (P1), (P2)
and $M(B)/B$ has
(P3).

It is clear that many such examples can be constructed.

Proposition \ref{IIILfulinf} is not needed in Example \ref{4E} but
 will be used again later.

\begin{Lem}\label{IIILu}
Let $A$ be a unital \CA\, and $0\le a\le 1$ be an element in $A.$
Suppose that there is $x\in A$ such that $x^*ax=1.$ Then there is
$N>0$ depends on $\|x\|$ (not on $A$ or $a$) for which there is
$y\in A$ with $\|y\|\le 1$ such that
$$
y^*f_N(a)y=1.
$$
In particular, $f_N(a)$ is full (where $f_N$ is as defined in
(\ref{fn})).
\end{Lem}

\begin{proof}
Let $q=a^{1/2}xx^*a^{1/2}.$ Then $q$ is a projection. There exists
$k>0$ depends on $\|x\|$ such that
$$
\|f_k(t)t^{1/2}-t^{1/2}\|<{1\over{16\|x\|^2}}\rforal t\in [0,1],
$$
where $f_k$ is as defined in (\ref{fn}).  Then
$$
\|f_k(a)q-q\|=\|(f_N(a)a^{1/2}-a^{1/2})x^*xa^{1/2}\|<1/16.
$$
It follows from A8 in \cite{Eff}  that there is a projection $p\in
\overline{f_k(a)Af_k(a)}$ such that
$$
\|q-p\|<1/2.
$$
Thus there exists $w\in A$ such that $w^*w=1$ and $ww^*=p.$ Choose
$N=k+1.$ Then $f_N(a)q=q.$ Thus
$$
w^*f_N(a)w=1.
$$
\end{proof}

\begin{Lem}\label{IILfulinf1}
Let $A$ be unital \CA\, and $0\le a\le 1$ be a full element in
$A.$ Suppose that there are $x_1,x_2,...,x_m\in A$ such that
$$
\sum_{i=1}^mx_i^*ax_i=1.
$$
Let $r=\sum_{i=1}^m\|x_i\|^2.$ Suppose also that
$1_{M_m(A)}\lesssim 1 .$ Then there exists an integer $N>0$
depends on $r$ (but not $A$ nor on $a$) such that $f_N(a)$ is
full. Moreover, there are $y_1,y_2,...,y_m\in A$ such that
$\sum_{i=1}^m\|y_i\|^2\le 1$ and
$$
\sum_{i=1}^m y_i^*f_N(a)y_i=1.
$$
\end{Lem}

\begin{proof}
Let
$$
X=\begin{pmatrix} x_1 & x_2 &\cdots & x_m\\
                  0 & 0 &\cdots & 0\\
                  \hdotsfor{4} \\
                  0 & 0 & \cdots & 0
                  \end{pmatrix} \andeqn
                  b=\begin{pmatrix} a & 0 &\cdots &0\\
                                    0 & a & \cdots & 0\\
                                    \hdotsfor{4} \\
                                    0& 0 &\cdots &a.
                                    \end{pmatrix}
                                    $$
Since $1_{M_m(A)}\lesssim 1,$ one obtain
                                    $Y\in M_m(A)$ with $\|Y\|=1$
                                    such that
$Y^*{\rm diag}(1,0,...,0)Y=1_{M_m(A)}.$   Note we have $0\le b\le
1$ and $XbX^*={\rm diag}(1,0,...,0).$ Thus
$$
Y^*XbX^*Y=1_{M_m(A)}.
$$
We compute that $\|X^*Y\|\le r^{1/2}.$ It follows from the above
lemma that there is $N>0$ $f_N(b)$ for which there is $z\in
M_m(A)$ with $\|z\|\le 1$ such that
$$
z^*f_N(b)z=1_{M_m(A)}.
$$
So $Yz^*f_N(b)zY^*=1.$ An easy computation shows that there are
$y_1,y_2,...,y_n\in A$ such that $\sum_{i=1}^m\|y_i\|^2\le 1$ and
$$
\sum_{i=1}^my_i^*f_N(a)y_i=1.
$$
\end{proof}

\begin{Prop}\label{IIILfulinf}
Let $\{A_n\}$ be a sequence of unital \CA s which has property
(P1). Then $l^{\infty}(\{A_n\})$ also has property (P1).
\end{Prop}

\begin{proof}
Let $a=\{a_n\}$ be a full element in $l^{\infty}(\{A_n\})$ such
that $0\le a\le 1.$ By \ref{IILfulinf1}, there exists $N>0$ such
that $f_N(a)$ is full. For each $n,$ there exists $x_n\in A_n$
such that $x_n^*f_N(a_n)x_n=1.$ Note that
$f_{N+1}(a_n)f_N(a)=f_N(a).$ It follows from \ref{PI2} that, for
each $n,$ there is $y_n\in A$ with $\|y_n\|\le 1$ such that
$$
y_n^*f_{N+1}(a)y_n=1.
$$
Put $y=\{y_n\}.$ Then $y\in l^{\infty}(\{A_n\}).$ It is clear that
there is $g\in C_0((0,1])_+$ such that
$$
\|g(a)ag(a)-f_{N+1}(a)\|<1/4.
$$
Then
$$
\|y^*g(a)ag(a)y-1\|=\|y^*(g(a)ag(a)-f_N(a))y\|\le 1/4.
$$
It follows that there is $z\in l^{\infty}(\{A_n\})$ with
$\|z\|<4/3$ such that
$$
z^*y^*g(a)ag(a)yz=1.
$$
\end{proof}

The above proposition is not required in the following example.
However it will be used in \ref{IILoinf}.

\begin{Ex}\label{4E}
{\rm Let $A$ be a unital separable  amenable purely infinite
simple \CA s. Denote by $B=c_0(A).$
Then $M(B)=l^{\infty}(A).$
 Put $q_{\infty}(A)=l^{\infty}(A)/c_0(A).$
So $M(B)/B=q_{\infty}(A).$

(1) $M(B)$ and $M(B)/B$ has property (P1) and (P2).

(2) $M(B)/B$ has property (P3).

It is clear that  (1) is obvious (it also follows from  \ref{IIILfulinf}).
 In fact, if $C=C_0((0,1), A),$ then $M(C)$ and $M(C)/C$
also have property (P1) and (P2). This could be proved rather
easily.

To see (2), let $D$ be a separable
\SCA\, of $M(B).$ Suppose that  $x^{(1)}=\{x_n^{(1)}\},$ $
x^{(2)}=\{x_n^{(2)}\},..., x^{(k)}=\{x_n^{(k)}\},...$ is
a dense sequence of the unit ball of $D.$
Using the fact that $A\otimes {\cal O}_{\infty}
\cong A$ (see Theorem 3.15 of \cite{KP}), we obtain a sequence
of \hm s $\phi_n: {\cal O}_{\infty}\to A$ such that
$$
\lim_{n\to\infty}\|\phi_n(b)a-a\phi_n(b)\|=0
$$
for all $a\in A$ and $b\in {\cal O}_{\infty}.$
Let $e_1\in {\cal O}_{\infty}$ be a proper projection.
There is an integer $n(1)>0$ such that
$$
\|\phi_{n(1)}(e_1)x_1^{(1)}-x_1^{(1)}\phi_{n(1)}(e_1)\|<1/2.
$$
There is a projection $e_2\in {\cal O}_{\infty}$ such that
$e_1e_2=e_2e_1=0$ and $1>e_1+e_2.$
There is $n(2)>0$ such that
$$
\|\phi_{n(2)}(e_j)x_l^{(i)}-x_l^{(i)}\phi_{n(2)}(e_j)\|<1/4,\,\,\,\,
i,j,l=1,2.
$$
Continuing in this fashion, we obtain a sequence of mutually
orthogonal nonzero projections $\{e_m\}\subset {\cal O}_{\infty}$
and a subsequence $\{n(m)\}$ such that
$$
\|\phi_{n(m)}(e_j)x_l^{(i)}-x_l^{(i)}\phi_{n(m)}(e_j)\|<1/2^m,\,\,\,\,\,\,
i,j,l=1,2,...,m.
$$
Put $p^{(j)}=\{\phi_{n(m)}(e_j)\}\in l^{\infty}(A),$ $j=1,2,.....$
Then $p_m^{(i)}p_m^{(i)}=0$ if $i\not=j.$
Moreover,
$$
\|\pi(p^{(j)})\pi(\{x^{(i)}\})-\pi(\{x^{(i)}\})\pi(p^{(j)})\|=0.
$$
This implies that
$$
\pi(p^{(j)})\pi(d)=\pi(d)\pi(p^{(j)}).
$$
Put $a_n^{(j)}=p^{(j)},$ $j=1,2,....$
This shows that $M(B)/B$ has property (P3).

It is clear, in fact, that $l^{\infty}(\{A_n\})/c_0(\{A_n\})$ has
property (P3 if each $A_n$ is a unital purely infinite simple \CA.
}
\end{Ex}

\section{ Amenable \morp s}



\begin{Lem}{\rm (cf. 2.3 of \cite{AAP}, see also 5.3.2 of \cite{Lnb})}\label{LI01}
Let $A$ be a separable \CA\, and $\psi: A\to {\mathbb C}$ be a
pure state. Denote also by $\psi$ the extension of $\psi$ on
${\tilde A}$ and put $L=\{a\in {\tilde A}: \psi(a^*a)=0\}.$ Then,
for any $\ep>0$ and any finite subset ${\cal F}\subset A,$ there
is $z_i\in {\tilde A}_+$ with $\|z_i\|=1 $  such that $z_i\not\in
L$,
\begin{equation}
z_{i+1}z_i=z_i,\,\,\,i=1,2,\andeqn
\|z_i(\phi(a)-a)z_i\|<\ep/2,\,\,\,i=1,2,3,
 \end{equation}
for all $a\in {\cal F}.$ Moreover, if $\{e_n\}$ is an approximate
identity for $A,$ then, for some large $N,$
\begin{eqnarray}
\|e_nz_ie_n(\phi(a)-a)e_nz_ie_n\|<\ep \andeqn e_nz_ie_n\not\in L
\end{eqnarray}
for all $a\in {\cal F}$ and all $n\ge N.$
\end{Lem}

\begin{proof}
To simplify notation, we may assume that ${\cal F}$ is a
subset of the unit ball of $A.$
Let
$$
N=\{a\in {\tilde A}: \phi(a)=0\}.
$$
Note that $L$ is a closed left ideal. Let $C$ be the hereditary
\SCA\, given by $L\cap L^*.$ As in the proof of 5.3.2 in
\cite{Lnb}, we have $z_1, z_2,z_3,\in {\tilde A}$ with $\|z_i\|=1$
($i=1,2,3$) such that $z_i\not\in L,$ $z_{i+1}y_i=z_i,$ $i=1,2,3,$
and
$$
\|z_i(\psi(a)-a)z_i\|<\ep/2,\,\,\,i=1,2,3.
$$
Let $\{e_n\}$ be an approximate identity for $A$ such that
$e_ne_{n+1}=e_n$ for all $n.$
Note $z_i$ has the form $\lambda_i 1_B+y_i',$ where $y_i'\in A$
and $\lambda_i\in {\mathbb C},$ $i=1,2.$
Choose large $n$ so that
$$
\|e_ka-ae_k\|<\ep/4,\,\,\, \|e_ka-a\|<\ep/4 \andeqn
\|e_kz_i-z_ie_k\|<\ep/4
$$
for all $a\in {\cal F}\cup \{z_1az_1, z_2az_2, z_3zz_3: a\in {\cal
F}\}$ and for all $k\ge n.$ Let $y_i=e_nz_ie_n.$ Then, for $n\ge
N,$
$$
\|y_i(\psi(a)-a)y_i\|<\ep/4+\|e_n^2z_i(\psi(a)-a)z_ie_n^2\|<\ep
\,\,\,{\rm for\,\,\,all}\,\,\, a\in {\cal F}.
$$
\end{proof}

The following is a folklore.

\begin{Lem}\label{LIapp}
Let $A$ be a \SCA\, of $B$ and $a\in A_+.$ Denote by $C$ the
hereditary \SCA\, of $B$ generated by $a.$ Then, for any
approximate identity $\{e_n\}$ of $A,$
$$
\|e_nb-b\|\to 0 \andeqn \|be_n-b\|\to 0,\,\,\,{\rm as}\,\,\, n\to\infty
$$
for all $b\in C.$
\end{Lem}

\begin{proof}
There exists a sequence of positive function $f_n\in C_0(sp(a))$ with
$0\le f_n\le 1$ such that $\{f_n(a)\}$ forms an approximate
identity for $C.$ Fix an element $b\in C.$ For any $\ep>0,$ there is $f_k$ such that
\begin{eqnarray}
\|f_k(a)b-b\|<\ep/4\andeqn \|bf_k(a)-b\|<\ep/4.
\end{eqnarray}
Choose integer $N>0$ such that
\begin{eqnarray}
\|e_nf_k(a)-f_k(a)\|<{\ep\over{4(\|b\|+1)}}\,\,\,{\rm for\,\,\,
all}\,\,\,n\ge N.
\end{eqnarray}
It follows that
\begin{eqnarray}
\|e_nb-b\| &\le &\|e_nb-e_nf_k(a)b\|+\|e_nf_k(a)b-f_k(a)b\|+\|f_k(a)b-b\|\\
  &<& \ep/4+\|b\|({\ep\over{4(\|b\|+1)}})+\ep/4=3\ep/4<\ep
\end{eqnarray}
\end{proof}

\begin{Lem}\label{LI1}

   Let $B$ be a unital \CA\, that has
the property (P1).
 Let $A$ be a  separable \CA\, and $I$ be an ideal of $A.$
Suppose that $j: A\to B$ is an  embedding such that
$j(a)$ is a full element of $B$ for all $a \not\in I.$
 Then, for any pure state $\phi: A\to {\mathbb C} 1_B\subset B$
which vanishes on $I,$
any finite subset ${\cal F}\subset A,$ and any $\ep>0,$
there is a partial isometry $V\in B$ such that
$$
\|\phi(a)-V^*j(a)V\|<\ep \,\,\,{\rm for}\,\,\,\,\,\, a\in {\cal F},\,\,\,
 V^*V=1_B \andeqn VV^*\in Her(j(A)).
$$
\end{Lem}

\begin{proof}
To simplify notation,  we identify $A$ with $j(A).$ Fix $0<\ep
<1/2.$
 By \ref{LI01}, there are
$z_1,z_2,z_3\in {\tilde A}_+$ with  $\|z_i\|=1$ and
$z_i\not\in L,$ $z_{i+1}z_i=z_i,$ $i=1,2$ such that
\begin{equation}\label{q1-01}
\|\phi(a)z_i^2-z_ij(a)z_i\|<\ep/4\,\,\,{\rm for}\,\,\, a\in {\cal F}
\,\,\,{\rm (}i=1,2{\rm )}
\end{equation}
for all $a\in {\cal F}.$
Note $L=\{a\in {\tilde A}: \psi(a^*a)=0\}.$
Therefore $I\subset L\cap L^*\subset L.$
Let $\{e_n\}$ be an approximate identity for $A$ such that
$e_ne_{n+1}=e_n,$ $n=1,2....$ Let $N$ be the integer as described in \ref{PI2}
so that
\begin{eqnarray}\label{e1-1}
\|\phi(a)(e_nz_ie_n)^2- e_nz_ie_nj(a)e_nz_ie_n\|<\ep/2,\,\,\,i=1,2,3.
\end{eqnarray}
Put $y_1=e_Nz_1e_N.$ We may assume that $y_1\not\in L.$
By the assumption, $y_1$ is full. Because $B$ has property (P1),  there
exists $x\in B$ such that $x^*y_1^2x=1_B.$ Put $v_1=y_1x.$ Then
$v_1^*v_1=1_B$ and $v_1v_1^*=p_1$ is a projection. Note that
$p_1\in Her(y_1).$ There is a projection in $q_1\in
Her(z_1^{1/2}e_Nz_1^{1/2})$ such that $q_1$ is equivalent to
$p_1.$ Therefore there is a partial isometry $w_1\in B$ such that
$w_1^*q_1w_1=1_B$ and $w_1w_1^*=q_1.$ Since $z_2^2z_1=z_1,$
$z_2^2q_1=q_1.$ By applying \ref{LIapp}, one can choose a large
integer $k>N$ so that, for all
 $n\ge k,$
\begin{eqnarray}\label{e1-02}
\|e_nq_1-q_1\|<\ep/32\andeqn \|e_nz_i-z_ie_n\|<\ep/32, \,\,\,i=1,2,3.
\end{eqnarray}
Thus
\begin{eqnarray}
\|(e_kz_2e_k)^2q_1-q_1\|=\|e_kz_2e_k^2z_2e_kq_1-q_1\|<8\ep/32=\ep/4.
\end{eqnarray}
Put $y_2=e_kz_2e_k.$ Then one estimates
\begin{eqnarray}
\|w_1^*y_2^2w_1-1\|=\|w_1^*q_1y_2^2q_1w_1-w_1^*q_1w_1\|<\ep/2.
\end{eqnarray}
Thus there is $s\in Her(z_1^{1/2}e_Nz_1^{1/2})_+\subset B_+$ such that $\|s\|\le {1\over{1-\ep/2}}$ and
\begin{eqnarray}
s^{1/2}w_1^*y_2^2w_1s^{1/2}=1.
\end{eqnarray}
Note that
\begin{eqnarray}\label{e1-032}
\|w_1s^{1/2}\|\le \sqrt{1\over{1-\ep/2}}<\sqrt{2\over{2-1/2}}=\sqrt{4/3}={2\sqrt{3}\over{3}}.
\end{eqnarray}
Define $V=y_2w_1s^{1/2}.$ Note that
\begin{eqnarray}\label{e1-03}
V^*V=1_B\andeqn VV^*\in Her(j(A)).
\end{eqnarray}
Put $y_3=e_{k+1}z_3e_{k+1}.$ Then, by (\ref{e1-02}),
\begin{eqnarray}\label{e1-0321}
\|y_3y_2-y_2\|=\|e_{k+1}(z_3e_kz_2-z_2)e_k\|<\ep/32.
\end{eqnarray}
Furthermore, by (\ref{e1-032}) and (\ref{e1-0321}),
\begin{eqnarray}\label{e1-04}
\|y_3V-V\|=\|y_3y_2w_1s^{1/2}-y_2w_1s^{1/2}\|\le (\ep/32)( {2\sqrt{3}\over{3}})=\sqrt{3}\ep/48.
\end{eqnarray}
We estimate, by applying (\ref{e1-04}),(\ref{e1-03}) and
(\ref{e1-1})
\begin{eqnarray}
\|\phi(a)-V^*aV\| &= &\|\phi(a)V^*V-V^*aV\| {\nonumber}\\
&\le& 3\sqrt{3}\ep/48+\|\phi(a)V^*y_2^2V-V^*y_2ay_2V\|{\nonumber}\\
&\le & 3\sqrt{3}\ep/48+\|\phi(a)y_2^2-y_2ay_2\| {\nonumber}\\
&<&3\sqrt{3}\ep/48+\ep/2<\ep
\end{eqnarray}
for all $a\in {\cal F}.$

\end{proof}

\begin{Remark}
{\rm  If $A$ has a unit, then the proof of Lemma \ref{LI1} is
almost identical to that of 5.3.2 of \cite{Lnb} which has its
origin in \cite{AAP}. When $A$ has no unit, elements $z_1, z_2,
z_3$ are not in $A_+$ but in ${\tilde{A}}_+.$ By using an
approximate identity $\{e_n\},$ one does have $\|y_3y_2-y_2\|$
small. However the norm $x$ could be large and depends on the
choice of $z_i$ as well as $N$ as in the above proof. By
introducing of $q_1,$ we are able to control the norm of
$w_1s^{1/2}.$}
\end{Remark}

\begin{Lem}\label{LI2}
Let $B$ be a unital \CA\, which has the property {\rm (}P1{\rm )}
and $A$ be a separable \CA. Suppose that there exists a sequence
of \hm\, $\phi_n: A\to B$ such that  $\{\phi_n(a): n=1,2,...\}$ is
a mutually orthogonal set in $B$ for all $a\in A.$ Let $I$ be an
ideal of $A$ such that ${\rm ker}\phi_n\subset I$ and $\phi_n(a)$
is a full element in $B$ for all $a\not\in I$ for all $n.$  Then,
for any state $\psi: A/I\to {\mathbb C}1_B\subset B,$ any finite
subset ${\cal F}\subset A,$ and any $\ep>0,$ there is a partial
isometry $V\in B$ and an integer $n$ such that
\begin{equation}
\|\psi\circ \pi(a)-V^*(\sum_{k=1}^n\phi_k(a))V\|<\ep\,\,\,\,
{\rm for\,\,\,} a\in {\cal F},\,\,\, V^*V=1B\andeqn VV^*\in Her(\sum_{i=1}^n\phi_i(A)),
\end{equation}
where $\pi: A\to A/I$ is the quotient map.
Moreover, if $\psi$ is only assume to be a nonzero positive linear
functional with $\|\psi\|\le 1,$ then the above still holds
where $V$ is merely a contraction.
\end{Lem}

\begin{proof}
By the Krein-Milman theorem, we have positive numbers
$\alpha_1,\alpha_2,...,\alpha_m$ with
$\sum_{i=1}^m \alpha_i=1$ and pure states $\psi_1,\psi_2,...,\psi_m$
of $A/I$
such that
\begin{equation}
\|\psi\circ \pi(a)-\sum_{i=1}^m\alpha_i\psi_i(a)\|<\ep/2\,\,\,\,{\rm for}\,\,\,
a\in {\cal F}.
\end{equation}
Let $\pi_n: A\to A/{\rm ker}\phi_n$
and $\gamma_n: A/{\rm ker}\phi_n\to A/I$ be the quotient maps, $n=1,2,....$
Note that $\psi_i\circ \gamma_n$
is a pure state of $A/{\rm ker}\phi_n.$

By \ref{LI2}, there are $V_i\in B$ such that
\begin{equation}
V_i^*V_i=1_B, V_iV_i^*\in Her(\phi_i(A))\andeqn
\|\psi_i\circ \gamma_i(\phi_i(a))-V_i^*\phi_i(a)V_i\|<\ep
\end{equation}
for all $a\in {\cal F}.$
One should note that
$$
\psi_i\circ \gamma_i\circ \phi_i=\psi_i\circ \pi.
$$

Set $V=\sum_{i=1}^m\sqrt{\alpha_i}\, V_i\in B.$
We see that $V^*V=\sum_{i=1}^m\alpha_i1_B=1_B$
and $VV^*=\sum_{i=1}^m\alpha_iV_iV_i^*\in Her(\sum_{i=1}^m\phi_i(A)).$
Moreover
\begin{eqnarray}
\|\psi\circ \pi(a)-
V^*(\sum_{i=1}^m\phi_i(a)V\| &=&
\|\psi\circ\pi(a)-\sum_{i=1}^m\alpha_iV_i^*\phi_i(a)V_i\|{\nonumber}\\
&\le & \|\psi\circ \pi(a)-\sum_{i=1}^m\alpha_i\psi_i(a)\|
+\sum_{i=1}^m\alpha_i\|\psi_i\circ \pi(a)-V_i^*\phi_i(a)V_i\|{\nonumber}\\
&<& \ep/2+\ep/2=\ep\,\,\,\,\,\,\,\,\,\,\,a\in {\cal F}.
\end{eqnarray}

To see the last statement of the lemma holds, we note that there
is $0<\lambda\le 1$ such that $\psi(a)=\lambda \cdot g(a)$ for
some state $g$ and for all $a\in A.$
\end{proof}

\begin{Lem}\label{LI3}
Let $A$ be a separable \CA\, and $B$ be a unital \CA\, which has
the property (P1) and (P2). Let $C$ be as described in \ref{Dp2}
with $n=k$ (see \ref{IIP2}).
Suppose that $\phi_n: A\to B$ is a sequence of \hm s such that
$\{\phi_n(a): n=1,2,...\}$ is a set of mutually orthogonal
elements in $B.$  Suppose that $I$ is an ideal of $A$ such that
$I\supset {\rm ker}\phi_n$ and $\phi_n(a)$ is a full element for
all $a\not\in I.$ Let $\psi: A/I\to M_k({\mathbb C})\subset
M_k(C)\subset B$ be a \morp\, Then for any finite subset ${\cal
F}\subset A$ and $\ep>0,$ there exists a contraction $V\in B$ and
an integer $m>0$ such that
\begin{eqnarray}
\|\psi(a)-V^*(\sum_{i=1}^m\phi_i(a))V\|<\ep \,\,\,\,{\rm for}\,\,\,
a\in {\cal F} \andeqn VV^*\in Her(\sum_{i=1}^K\phi_i(A)),
\end{eqnarray}
where $\pi: A\to A/I$ is the quotient map.
\end{Lem}

\begin{proof}
Write $\psi(a)=\sum_{i=1}^k\psi_{ij}(a)\otimes s_{ij}$ for
$a\in A,$ where $\{s_{ij}\}$ is a system of matrix units for
$M_n$ and $\psi_{ij}: A\to {\mathbb C}$ is linear.
Note we also assume that $s_{ii}$ are as in \ref{Dp2} and \ref{IIP2}, $i=1,2,...,k.$
Define $\Phi: M_k(A)\to {\mathbb C}\subset C$ by
$\Phi((a_{ij})_{k\times k})=\sum_{i,j=1}\psi_{ij}(a_{ij}),$ where
$a_{ij}\in A.$
Let $Z$ be as in \ref{Dp2} so that $ZZ^*=s_{11}.$
Put $J_n(a)=Z_k\phi_n(a)Z^*$ for all $a\in A.$ So
$J_n$ maps $A$ into $C$ ($=s_{11}Bs_{11}$).
Note that  $\phi_n\otimes {\rm id}: M_k(A)\to
M_k(C)$ is also full.
Set ${\cal G}=\{(a_{ij}): a_{ij}\in {\cal F}\cup \{0\}\}.$
Thus, by applying \ref{LI2}, there is $W\in M_k(B)$ with
$\|W\|\le 1$ such that
\begin{eqnarray}
\|\Phi(b)-W^*(\sum_{k=1}^mJ_k\otimes {\rm id}(b))W\|<\ep/2n^2
\,\,\,{\rm for}\,\,\,
b\in {\cal G}.
\end{eqnarray}
Note that we may also assume that
$W^*W\le {\rm diag}(1_C,0,...,0).$
Choose a positive element
$0\le d\le 1$ in $A$ such that
\begin{eqnarray}
\|da-a\|<\ep/2n^2\,\,\,{\rm for\,\,\, all}\,\,\, a\in {\cal F}.
\end{eqnarray}
Let $v_i=(\overbrace{0,...,0}^{i-1},d,0,...,0)$ and
$v_i'$ be the $n\times n$ matrix with the first row as $v_i$ and
rest row are zero.
Put $r_i=\sum_{n=1}^mJ_n\otimes {\rm id}(v_i').$ Note that, for any
$a\in A,$
\begin{eqnarray}
\|r_i^*(\sum_{n=1}^mJ_n(a)\otimes {\rm id}(a\otimes s_{11}))r_j-
\sum_{n=1}^mJ_n\otimes {\rm id} (a\otimes s_{ij})\|<\ep/2n^2.
\end{eqnarray}
Therefore
\begin{eqnarray}
\|\psi_{ij}(a)- W^*r_i^*(\sum_{k=1}^mJ_k\otimes {\rm id}(a\otimes
s_{11})) r_jW\|<\ep/n^2
\end{eqnarray}
for all $a\in {\cal F}.$
Put $V'=(v_1'W,v_2'W,...,v_n'W).$ Note we view $V'$ is an $n\times
n$ matrix with $i$-th column as a nonzero column of $v_i'W,$
$i=1,2,...,n.$ Then
\begin{eqnarray}
\|\psi(a)-V'^*\sum_{k=1}^mJ_k\otimes {\rm id}(a\otimes e_{11})V'\|<\ep
\,\,\,{\rm for}\,\,\, a\in {\cal F},
\end{eqnarray}
Define $V=Z^*V',$ we have
\begin{eqnarray}
\|\psi(a)-V^*\sum_{n=1}^m\phi_n(a)V\|<\ep \,\,\,{\rm for}\,\,\, a\in {\cal F}.
\end{eqnarray}
We also note that
$VV^*\in Her(\sum_{n=1}^m\phi_n(A)).$

\end{proof}

\begin{Lem}\label{LIf}
Let $A$ be a separable \CA\, and $B$ be a unital \CA\, which has
the property (P1) and property (P2). Suppose that $\phi_n: A\to B$
is a sequence of \hm s such that the embedding $j_n: \phi_n(A)\to
B$  is full where $\{\phi_n(a): n=1,2,...\}$ is a set of mutually
orthogonal elements in $B.$ Suppose that $\psi: A\to B$ is
amenable such that ${\rm ker}\psi\supset {\rm ker}\phi_n,$
$n=1,2,....$  Then, for any finite subset ${\cal F}\subset A$ and
$\ep>0,$ there exists a contraction $V\in B$ and an integer $K>0$
such that
\begin{eqnarray}
\|\psi(a)-V^*(\sum_{i=1}^K\phi_i(a))V\|<\ep \,\,\,\,{\rm for}\,\,\,
a\in {\cal F} \andeqn VV^*\in Her(\sum_{i=1}^K\phi_i(A)).
\end{eqnarray}
\end{Lem}

\begin{proof}
Fix a finite subset ${\cal F}$ and $\ep>0.$ Since $\psi$ is
amenable, to simplify notation, without loss of generality, we may
assume that $\psi=\alpha\circ \beta,$ where $\beta: A\to
M_n=M_n({\mathbb C}\cdot 1_C)$ and $\alpha: M_n\to B$ are \morp s
(it should be noted though that $n$ depends on ${\cal F}$ as well
 as $\ep$).
Write $M_n(C)\subset B$ as in \ref{Dp2} (see also \ref{IIP2}).
Put ${\cal G}=\beta({\cal
F}).$ It is convenient to assume that ${\cal F}$ lies in the unit
ball of $A$ so ${\cal G}$ lies the unit ball of $M_n({\mathbb
C}\cdot 1_C).$ Note that $\sigma: M_n\to M_n({\mathbb C})\subset
B$ is full. There exists an integer $m>0$ and a contraction $Z\in
M_m(B)$ such that
\begin{eqnarray}
\|\alpha(b)-Z^*{\rm
diag}\overbrace{(b,b,...,b)}^{m}Z\|<\ep/4\,\,\,{\rm for}\,\,\,
b\in {\cal G}.
\end{eqnarray}
It follows from \ref{LI3} that
there is $N(1)>1$ and a contraction $W_1\in B$ such that
\begin{eqnarray}
\|\beta(a)-W_1^*\sum_{i=1}^{N(1)}\phi_i(a)W_1\|<\ep/4m\,\,\,\,
{\rm for}\,\,\, a\in {\cal F}
\end{eqnarray}
as well as integers $N(k+1)>N(k)$ and a contractions $W_k\in B$ such that
\begin{eqnarray}
\|\beta(a)-W_{k+1}^*\sum_{i=N(k)+1}^{N(k+1)} \phi_i(a)W_{k+1}\|<\ep/4m
\,\,\,\,
{\rm for}\,\,\, a\in {\cal F}, k=1,2,...,.
\end{eqnarray}
Note
we have
\begin{eqnarray}
\|\alpha\circ \beta(a)-Z^*{\rm
diag}\overbrace{(\beta(a),\beta(a),...,\beta(a))}^{m}Z\|<\ep/2\,\,\,{\rm
for}\,\,\, a\in {\cal F}.
\end{eqnarray}
It follows that
$$
\|\psi(a)-Z^*({\rm diag}(W_1^*\sum_{i=1}^{m(1)}\phi_i(a)W_1,
\cdots, W_m\sum_{i=N(m-1)+1}^{N(m)}\phi_i(a))W_m))Z\|<\ep/2
$$
for all $a\in {\cal F}.$
There exists $d_i\in Her(\phi_i(A))_+$ with $0\le d_i\le 1$
such that
\begin{eqnarray}
\|d_i\phi_i(a)-\phi_i(a)\|<\ep/2m\andeqn \|d_i\phi_i(a)d_i-\phi_i(a)\|
<\ep/2m
\end{eqnarray}
for all $a\in {\cal F}.$
Note that $d_id_j=0$ if $i\not=j,$ $i,j=1,2,...,m.$
Now let $Y$ be the $n\times n$ matrix
so that the first row is $(d_1, d_2,...,d_m)$ and the rest are zero.
Put $W={\rm diag}(W_1,W_2,...,W_m)$ and
$V=YWZ.$
Then
$$
\|{\rm diag}(W_1^*\sum_{i=1}^{m(1)}\phi_i(a)W_1, \cdots,
W_m\sum_{i=N(m-1)+1}^{N(m)}\phi_i(a)W_m)
-W^*Y^*\sum_{k=1}^{N(m)}\phi_k(a)YW\|<\ep/2
$$
for $a\in {\cal F}.$ Moreover
\begin{eqnarray}
\|\psi(a)-V^*\sum_{k=1}^{N(m)}\phi_k(a)V\|<\ep \,\,\,{\rm
for}\,\,\,a\in {\cal F}\andeqn VV^*\in
Her(\sum_{k=1}^{N(m)}\phi_k(A)).
\end{eqnarray}
\end{proof}

\section{Commutants in the ultrapower of corona algebras}

\begin{Def}\label{IIDfilter}
{\rm Recall that a family ${\omega}$  of subsets of ${\mathbb N}$
is an {\it ultrafilter} if

(i) $X_1,...,X_n\in \omega$ implies $\cap_{i=1}^n X_i\in \omega,$

(ii) $\O\not\in \omega,$

(iii) if $X\in\omega$ and $X\subset Y,$ then $Y\in\omega$ and

(iv) if $X\subset {\mathbb N}$ then either $X$ or ${\mathbb
N}\setminus X$ is in $\omega.$

An ultrafilter is said to be {\it free}, if $\cap_{X\in \omega} X
=\emptyset.$ The set of free ultrafilters is identified with elements in
$\beta {\mathbb N}\setminus {\mathbb N},$ where $\beta{\mathbb N}$
is the Stone-Cech compactification of ${\mathbb N}.$

 A sequence $\{x_n\}$  (in a normed space ) is
said to converge to $x_0$ along $\omega,$ written
$\lim_{\omega}x_n=x_0,$ if for any $\ep>0$ there exists
$X\in\omega$ such that $\|x_n-x_0\|<\ep$ for all $n\in X.$

Let $\{B_n\}$ be a sequence of \CA s.
Fix an ultrafilter $\omega.$ The ideal of $l^{\infty}(\{B_n\})$
which consists of those sequences $\{a_n\}$ in $l^{\infty}(\{B_n\})$
such that $\lim_{\omega}\|a_n\|=0$ is denoted by $c_{\omega}(\{B_n\}).$
Define
$$
q_{\omega}(\{A_n\})=l^{\infty}(\{B_n\})/c_{\omega}(\{B_n\}).
$$
If $B_n=B, n=1,2,...,$ we use
$c_{\omega}(B)$ for $c_{\omega}(\{B_n\})$
and $q_{\omega}(A)$ for $q_{\omega}(\{A_n\}),$ respectively. }
\end{Def}

\begin{Lem}\label{IIL2}
Let $A$ be a \CA, $I$ be an ideal of $A$ and let
$a\in A\setminus \{0\}$ such that $0\le a\le 1.$
Suppose that $a\not\in I.$ Then there is $b\in C^*(a)$
with  $0\le b\le 1$ and $\|b\|=1$ such that if $c\in C^*(b)\setminus J,$ then $c\not
\in I,$ where
$$
J=\{f(b): f\in C_0(sp(b)\setminus \{0\}),\,\,\, f(1)=0\}.
$$
\end{Lem}

\begin{proof}
Let $\pi: A\to A/I$ be the quotient map.
Then $\pi(a)\not=0.$ Suppose that $\xi\in sp(a)\setminus \{0\}.$
Let $f\in C_0(sp(a)\setminus \{0\}) $ such that $f(\xi)=1$ and
$0<f(t)<1$ for all other $t\in sp(a).$
Set $b=f(a).$ Then, $\pi(b)\not=0.$ and $\|\pi(b)\|=1.$
If $c\not\in J,$ $c=g(a)$ for some $g\in C_0(sp(a)\setminus \{0\})$
such that $g(\xi)\not=0.$  It follows that $\pi(g(a))\not=0.$
Therefore $c\not\in I.$
\end{proof}

\begin{Lem}\label{IIL2r}
Let $B$ be a unital \CA\, and $a\in B$ be an element with $0\le
a\le 1.$ Suppose that there is $x\in B$ such that $x^*ax=1.$ Then
there exists an element $b\in C^*(a)$ such that $c$ is full for
all $c \in C^*(b)\setminus J,$ where
$$
J=\{f(b): f\in C_0(sp(b)\setminus \{0\}),\,\,\, f(1)=0\}.
$$
\end{Lem}

\begin{proof}
Put $v=a^{1/2}x.$ Then $v^*v=1$ and $vv^*=q$ for some projection $q\in B.$
Note that $q\in Her(a^{1/2}xx^*a^{1/2})\subset Her(a).$
For any $0<\ep/<1/4,$ there is $N>0$ such that
$$
\|f_n(a)p-p\|<\ep/2\,\,\,{\rm for\,\,\,all}\,\,\, n\ge N.
$$
($f_n$ be as in (\ref{fn}).) It follows that
$$
\|f_n(a)pf_n(a)-p\|<\ep
$$
for all $n\ge N.$ If follows that there is a projection $q\in
Her(f_N(a))$ and partial isometry $w\in B$ such that $w^*qw=1$ and
$ww^*=q.$ Thus $f_{N+1}(a)q=q.$ Put $b=f_{N+1}(a).$ Thus, for any
function $g\in C_0((0,1]),$ if $g(1)\not=0,$ then $g(b)q=q.$ It
follows that $w^*g(b)w=1.$ Thus $g(b)$ is full. The lemma follows.
\end{proof}

\begin{Lem}\label{IIILnor}
Let $A$ be a unital separable \SCA\, of a unital \CA\, $B$ which
has property (P1) and (P3). Suppose that every nonzero element in
$A$ is full in $B.$ Then there exists a sequence of sequences of
positive elements $\{a_n^{(i)}\},$ $i=1,2,...$ with $0\le a\le 1$
satisfying the following:

{\rm (1)} $\lim_{n\to\infty}\|a_n^{(i)}a-aa_n^{(i)}\|=0$ for all
$a\in A$ and $i=1,2,...;$

{\rm (2)} $\lim_{n\to\infty}\|a_n^{(i)}a_n^{(j)}\|=0$ if $i\not=j$
and

{\rm (3)} $\Pi(\{a_n^{(i)}\})\Pi\circ J(a)$ is full in
$q_{\omega}(A)$ for any free ultrafilter $\omega\in \beta{\mathbb
N}\setminus {\mathbb N},$ where $J: B\to l^{\infty}(B)$ is defined
by $J(b)=(b,b,...,b,...)$ for $b\in B$ and $\Pi: l^{\infty}(B)\to
q_{\omega}(B)$ is the quotient map.
\end{Lem}

\begin{proof}
For each nonzero element $0\le a\le 1$ in $A,$ define
$$
r(a)=\inf\{ \|x\|: x^*ax=1\}.
$$
 Let $b_1,b_2,...,b_n,...$ be a dense sequence of the
unit ball of $A.$ We may assume that $\{b_n\}$ contains a
subsequence of positive elements which is dense in the positive
part of the unit ball. For each $0\le b_k\le 1$ in the sequence,
from the assumption, there is $x_k\in B$ such that $x_k^*b_kx_k=1$
and $\|x_k\|\le (4/3)r(b_k).$  Let $D$ be the separable \SCA\,
generated by $A$ and $\{x_k\}.$

We claim that, for each nonzero $a\in A$ with $0\le a\le 1$ there
is $x\in D$ such that $x^*ax=1.$ There is $z\in B$ such that
$z^*az=1$ and $\|z\|<(4/3)r(a).$ There is $b_k$ with $0\le b_k\le
1$ for which
$$
\|a-b_k\|<1/8((4/3)r(a)+1)^2.
$$
Then
$$
\|z^*b_kz-1\|\le \|z^*(b_k-a)z\|\le 1/8.
$$
We obtain $y\in D$ with $\|y\|<8/7$ such that
$$
y^*z^*b_kzy=1.
$$
It follows that $r(b_k)\le (8/7)r(a).$ Hence there is $x_k\in D$
with $\|x_k\|\le (4/3)(8/7)r(a)$ such that $x_k^*b_kx_k=1.$ It
follows that
$$
\|x_k^*ax_k-1\|\le
\|x_k^*(a-b_k)x_k\|<(1/8((4/3)r(a)+1)^2)[(4/3)(8/7)r(a)]^2<8/49<1.
$$
Thus there is $d\in D$ such that
$$
d^*x_kax_kd=1.
$$
This proves the claim.

Now since $B$ has property (P3) and $D$ is separable, there exists
a sequence of sequences of nonzero elements $\{a_n^{(i)}\}$ in $B$
with $0\le a_n^{(i)}\le 1$ such that

(i) $\lim_{n\to\infty}\|a_n^{(i)}d-da_n^{(i)}\|=0$ for all $d\in
D;$

(ii) $\lim_{n\to\infty}\|a_n^{(i)}a_n^{(j)}\|=0$ if $i\not=j$ and

(iii) for each $i,$ $\{a_n^{(i)}\}$ is full in $l^{\infty}(A).$

Thus (1) and (2) follow. To see (3), let $a\in A.$ From the claim,
there is $d\in D$ such that
$$
d^*ad=1.
$$
Put $a_i=\{a_n^{(i)}\}.$ Then, by \ref{IIILfulinf}, there is $z\in
l^{\infty}(A)$ such that $z^*a_iz=1.$  Note (i) implies that
$$
\Pi(a_i)\Pi\circ J(b)=\Pi\circ J(b)\Pi(a_i)\rforal b\in D.
$$
Put $g=\Pi\circ J(d)\Pi(z).$ Then
$$
g^*\Pi(a_i)\Pi\circ J(a)g=\Pi(z^*)\Pi\circ J(d^*)\Pi(a_i)\Pi\circ
J(a)\Pi\circ J(d)\Pi(z)
$$
$$
=\Pi(z^*)\Pi(a_i)\Pi\circ J(d^*)\pi\circ J(a)\Pi\circ J(d)\Pi(z)
=\pi(z^*)\Pi(a_i)\Pi(z)=1.
$$
\end{proof}

\begin{Lem}\label{IILoinf}
Let $A$ be a unital separable amenable \CA, $B$ be a unital \CA\,
which has  property (P1), (P2) and (P3). Let $\omega\in \beta
{\mathbb N}\setminus {\mathbb N}$ be a free ultrafilter. Suppose
that $\tau: A\to B$ is a full unital embedding. Let
$\tau_{\infty}: A\to l^{\infty}(B)$ be defined by
$\tau_{\infty}(a)=(\tau(a),\tau(a),...)$ and let $\psi=\Pi\circ
\tau_{\infty},$ where $\Pi: l^{\infty}(B) \to q_{\omega}(B).$ Then
there is a unital \SCA\, $C\cong {\cal O}_{\infty}$ in the
commutant of $\psi(A)$ in $q_{\omega}(B).$
\end{Lem}

\begin{proof}
Let $\{a_n^{(i)}\}$ be the sequence of sequences of elements given
by \ref{IIILnor}. Put $a_i=\{a_n^{(i)}\},$ $i=1,2,....$ Let $D$ be
as in the proof \ref{IIILnor}.

Applying \ref{IIL2r} (and also using $D$ as in the proof of
\ref{IIILnor}), we may assume that each $a_i$ has the property
that $sp(a_i)\subset [0,1]$ and $f(\Pi(a_i))$ is full for all
$0\le f\le 1$ in $C_0((0,1])$ for which $f(1)\not=0.$

Let $X=(0,1]$ and fix $i.$ Define $\phi_j', L': C_0(X)\otimes A\to
q_{\omega}(B)$ by $\phi_i'(f\otimes a)=f(\Pi(a_i))\psi(a)$ and
$L'(f\otimes a)=f(1)\psi(a)$ for $a\in A,$ respectively.  By (3)
in \ref{IIILnor}, $\phi_i'$ is full.  Let $\{{\cal F}_j\}$ be an
increasing sequence of finite subsets of $A$ for which
$\cup_{n=1}^{\infty} {\cal F}_j$ is dense in $A$ and $\{g_n\}$ be
a dense sequence of $C_0((0,1]).$

Let $\{a_{i(k)}\}_{k=1}^{\infty}$ be a subsequence of $\{a_i\}.$
It follows from \ref{LIf} that there exists $s_n\in B$ such that
\begin{eqnarray}
\|s_n^*(\sum_{k=1}^{m(n)}g_j(a_{i(k)})\psi(a))
s_n-g_j(1)\psi(a)\|<1/2^n \,\,\,{\rm for}\,\,\, a\in {\cal
F}_n\andeqn j=1,2,...,n
\end{eqnarray}
Moreover, $s_ns_n^*\in {\rm
Her}(\sum_{k=1}^{m(n)}(a_{i(k)})\psi(A)).$ Suppose that
$s_n=\Pi((s_{n,1},s_{n,2},...)),$ $n=1,2,....$ We may assume that
$$
\|s_{n,k(n)}^*(\sum_{k=1}^{m(n)}g_j(a_{k(n)}^{(i(k)})\tau(a))s_{n,k(n)}
-g_j(1)\tau(a)\|<1/2^n, n=1,2,....
$$
Now put $t_n=s_{n,k(n)},$ $t'=(t_1,t_2,...)$ and $t=\Pi(t').$
Define $\Phi: C_0(X)\otimes A\to l^{\infty}(B)$ by
$$
\Phi(f\otimes
a)=\{\sum_{k=1}^{m(n)}f(a_{k(n)}^{(i(k)})\tau(a)\}\rforal f\in
C_0(X)\andeqn a\in A.
$$
It follows that
$$
t^*\Pi\circ \Phi(f\otimes a)t=f(1)\psi(a)\rforal f\in
C_0(X)\andeqn a\in A.
$$
Put $b(\{i(k)\})=\Pi(\{a_{k(n)}^{i(k)}\}).$ Note that $0\le
b(\{i(k)\})\le 1.$ We have (with $\imath(t)=t$ for all $t\in
(0,1]$)
$$
t^*b(\{i(k)\})t=\imath(1)=1_{q_{\omega(B)}}.
$$
Put $w(\{i(k)\})=b(\{i(k)\})^{1/2}t$ and
$q=b(\{i(k)\})^{1/2}tt^*b(\{(i(k)\})^{1/2}.$ Since $b(\{i(k)\})\in
\psi(A)'$ and $\imath(1)=1,$ we have
\begin{eqnarray}\label{eoinf}
t^*b(\{i(k)\})^{1/2}\psi(a)b(\{i(k)\})^{1/2}t=t^*b(\{i(k)\})\psi(a)t=\imath(1)\psi(a)
=\psi(a)\rforal a\in A.
\end{eqnarray}
It follows from 6.36 in \cite{Rb} that
$w(\{i(k)\})=b(\{i(k)\})^{1/2}t\in \psi(A)'.$  It clear that if
$\{i(k)\}$ and $\{i(k)'\}$ are two disjoint infinite subsets of
${\mathbb N},$ then corresponding projections $q$ and $q'$ are
orthogonal. This implies that one has a sequence of isometries
$v_k\in \psi(A)'$ such that $v_k^*v_k=1_{q_{\omega}(B)}$ and $1\ge
\sum_{k=1}^nv_kv_k^*,$ $n=1,2,....$ Thus $\psi(A)'$ admits a
unital embedding of ${\cal O}_{\infty},$
\end{proof}

\section{Full extensions}

\begin{Def}\label{DKL}

{\rm
 Let ${\bf Ext}(A,B)$ be the usual set of {\it stable}
unitary equivalence classes of extensions of the form
(\ref{IVext}). When $A$ is amenable, it is known
(Arveson/Choi-Effros) that ${\bf Ext}(A,B)$ is a group. Moreover,
it can be identified with $KK^1(A, B).$ Let ${\cal T}(A,B)$ be the
set of all {\it stable} unitary equivalence classes of
approximately trivial extensions. It is known that ${\cal T}(A,B)$
is a subgroup of $KK^1(A,B)$ (see \cite{Lnuct}). Following R\o
rdam, one defines $KL^1(A,B)=KK^1(A,B)/{\cal T}(A,B).$

Let $G_i,$ $i=1,2,3$ be three abelian groups. A group extension
$0\to G_1\to G_3\to G_2\to 0$ is said to be {\it pure} if
every finitely generated subgroup of $G_2$ lifts.
Denote by $Pext(G_2, G_1)$ the set of all pure extensions and
$E(G_2, G_1)=ext_{\mathbb Z}(G_2, G_1)/Pext(G_2, G_1).$

If $A$ satisfies the Approximate Universal Coefficient Theorem
(AUCT) --see \cite{Lnuct}, then one has the following short exact
sequence:
\begin{eqnarray}\label{IVuct}
0\to E(K_i(A), K_i(B))\to KL^1(A, B)\to Hom(K_i(A),K_{i-1}(B))\to
0.
\end{eqnarray}
So $KL^1(A,B)$ is computable in theory. It should be noted every
separable amenable \CA\, which satisfies the Universal Coefficient
Theorem (UCT) satisfies the AUCT. Rosenberg and Schochet
(\cite{RS}) show that every separable \CA s in the so-called
``bootstrap" class satisfies the UCT (therefore the AUCT). We also
use the notation $KL(A,B)=KL^1(A,SB).$

As mentioned in the introduction, two stably unitarily equivalent
extensions are in general not unitarily equivalent and trivial
extensions are not unitarily equivalent. Furthermore, an essential
extension which is  zero in $KK^1(A, B)$ may not be trivial (or
approximately trivial). We will use $KL^1(A, M(B)/B)$ to give a
classification of full essential extensions up to approximately
unitary equivalence. }

\end{Def}

\begin{Prop}\label{IIPo2}
Let $D$ be a unital \CA\, for which there is a unital embedding
from ${\cal O}_2$ to $D.$ Let $h_1, h_2: {\cal O}_2\to D$ be two
full \hm s. Suppose that $h_1(1_{{\cal O}_2})\sim h_2(1_{{\cal
O}_2}).$ Then there is a sequence of partial isometries $v_n$ such
that
\begin{eqnarray}
v_n^*v_n=h_2(1_{{\cal O}_2}),\,\,\, v_nv_n^*=h_1(1_{{\cal
O}_2})\andeqn \lim_{n\to\infty}\|v_n^*h_1(a)v_n-h_2(a)\|=0
\end{eqnarray}
for all $a\in {\cal O}_2.$
\end{Prop}

\begin{proof}
This is the combination of Theorem 6.5 and Lemma 7.2 in
\cite{Lnsemi}.
\end{proof}


\begin{Lem}\label{IIILao}
Let $A$ be a unital separable \CA,  $B$ and $C$ be  unital \CA s
such that $B\otimes {\cal O}_2$ is a unital \SCA\, of $C$ and $C$
has property (P1). Suppose that $h_1, h_2: A\to B\otimes {\mathbb
C}\cdot 1\subset B\otimes {\cal O}_2$ are two unital full
monomorphisms. Then $h_1$ and $h_2$ are approximately unitarily
equivalent in $C.$
\end{Lem}

\begin{proof}
It follows from \cite{R3} that ${\cal O}_2\cong {\cal O}_2\otimes
{\cal O}_2 .$ Let $p_n=1_B\otimes q_n\otimes 1_{{\cal O}_2},$
where $\{q_n\}$ is a sequence of mutually orthogonal non-zero
projections in ${\cal O}_2.$ Note that $p_n\sim 1_{B\otimes{\cal
O}_2\otimes {\cal O}_2},$ $n=1,2,....$ Define
$\phi_i(a)=p_ih_1(a)$ and $\psi_i(a)=p_ih_2(a)$ for all $a\in A.$
Also define $\Phi_n(a)=(1-\sum_{i=1}^np_i)h_1(a)$ and
$\Psi_n(a)=(1-\sum_{i=1}^np_i)h_2(a)$ for all $a\in A.$ Then, for
each $n,$ $h_1=\sum_{i=1}^n \phi_i\oplus \Phi_n$ and
$h_2=\sum_{i=1}^n\psi_i\oplus \Psi_n.$
 Note that $\phi_i, \Phi_n,
\psi_i$ and $\Psi_n$ are all full.
Now we work in $B\otimes {\cal O}_2\otimes 1.$
There are partial isometries $v_{i,j}\in {\cal O}_2$ such that
\begin{eqnarray}
v_{i,j}^*v_{i,j}=p_j\andeqn v_{i,j}v_{i,j}^*=p_i, \,\,\,i,j=1,2,...,n\\
\andeqn v_{n+1, j}^*v_{n+1, j}=p_j,
\,v_{n+1,j}v_{n+1,j}^*=1-\sum_{i=1}^np_i, \,\,\,j=1,2,...,n.
\end{eqnarray}
Put $w_{i,j}=1\otimes v_{i,j}\otimes 1.$ Then we also have
\begin{eqnarray}
w_{i,1}^*\phi_1w_{i,1}=\phi_i,\,\,\, i=1,2,...,n
\andeqn w_{n+1,1}^*\phi_1w_{n+1,1}=\Phi_n.
\end{eqnarray}

Let ${\cal F}_1, {\cal F}_2,...,{\cal F}_n,...$ be an increasing
sequence of finite subsets of $A$ such that $\cup_{n=1}^{\infty}
{\cal F}_n$ is dense in $A.$ It follows from Lemma 5.4.2 of
\cite{Lnb} that, for each $n,$ there are isometries $u_n,v_n\in
B\otimes {\cal O}_2\otimes 1$ such that
$$
\|u_n^*h_1(a)u_n-h_2(a)\|<1/n\andeqn
\|v_n^*h_2(a)v_n-h_1(a)\|<1/n\,\,\,\,{\rm for}\,\,\, a\in {\cal
F}_n.
$$
 Note that the relative commutant of $B\otimes {\cal
O}_2\otimes 1$ contains a unital \SCA\, $1_B\otimes 1_{{\cal
O}_2}\otimes {\cal O}_2$ which is isomorphic to ${\cal O}_2.$ It
follows from 1.10 in \cite{KP} that $h_1$ and $h_2$ are
approximately unitarily equivalent.
\end{proof}

\begin{Lem}\label{IILo22}
Let $A$ be a unital separable nuclear \CA, $B_1, B_2$ be two
unital \CA\, and $C$ be another unital \CA. Suppose that $j_i:
B_i\otimes {\cal O}_2\to C$ are two full monomorphisms so that
$j_1(1)\sim j_2(1)$ and $h_i: A\to B_i$ are two full unital
monomorphisms. Then there is a sequence of partial isometries
$v_n\in C$ such that $v_n^*v_n=j_1(1),$ $v_nv_n^*=j_2(1)$ and
\begin{eqnarray}
\lim_{n\to\infty}\|v_n^*(j_2\circ h_2(a))v_n-j_1\circ h_1(a)\|=0\,\,\,
{\rm for\,\,\, all}\,\,\, a\in A.
\end{eqnarray}
\end{Lem}

\begin{proof}
To simplify notation, we may assume that $j_1(1)=j_2(1).$
Therefore we may assume that both $j_1$ and $j_2$ are unital.
Define $J_i: B\otimes {\cal O}_2\to l^{\infty}(C)$ by
$J_i(b)=(j_i(b), j_i(b),...)$ for $b\in B_i\otimes {\cal O}_2$ and
$H_i=J_i\circ h_i,$ respectively, $i=1,2.$ Note that  these maps
are full in $l^{\infty}(C).$ Since there is a unital ${\cal O}_2$
embedding to $l^{\infty}(C),$ by \ref{IIPo2}, we obtain unitaries
$u_n\in C$ such that
\begin{eqnarray}
\lim_{n\to\infty}\|u_n^*J_2(1\otimes b)u_n-J_1(1\otimes b)\|=0\,\,\,
{\rm for\,\,\, all}\,\,\, b\in {\cal O}_2.
\end{eqnarray}
Denote $U=\{u_n\}$ in $l^{\infty}(C).$ Let $\omega$ be a free
ultrafilter on ${\mathbb N}$ and $\pi: l^{\infty}(C)\to
q_{\omega}(C)$ be the quotient map. Let  $D$ be the \SCA\,
generated by $\pi\circ J_1(B_1\otimes {\mathbb C}\cdot 1_{{\cal
O}_2})$ and $\pi\circ {\rm ad}\,U\circ J_2(B_2\otimes {\mathbb
C}\cdot 1_{{\cal O}_2}).$ It follows that  $D',$ the commutant of
$D,$ contains $J_1(1_{B_1}\otimes {\cal O}_2)$ which is isomorphic
to ${\cal O}_2.$ Therefore we may write $D\subset D\otimes {\cal
O}_2.$ Now $\pi\circ H_1$ and $\pi\circ {\rm ad} W\circ H_2$ are
two full unital monomorphisms from $A$ into $D\subset D\otimes
{\cal O}_2.$ It follows from \ref{IIILao} that $\pi\circ H_1$ and
$\pi\circ {\rm ad}W \circ H_2$ are approximately unitarily
equivalent. It follows from Lemma 6.2.5 of \cite{Rb}  that
$j_1\circ h_1$ and $j_2\circ h_2$ are approximately unitarily
equivalent.
\end{proof}

\begin{thm}\label{IIITab}
Let $A$ be a unital separable nuclear \CA, $B$ be a unital \CA\,
which has property (P1), (P2) and (P3). Let $j_o: A\to {\cal
O}_2\to B$ be a full embedding of $A$ into $B$ which factors
through ${\cal O}_2.$ Suppose that $\tau: A\to B$ is a full
monomorphism. Then there is a sequence of partial isometries
$V_n\in M_2(B)$ such that $V_n^*V=1_B\oplus j_o(1_A),$ $
V_nV_n^*=1_{B}$ and
$$
\lim_{n\to\infty}\|V_n(\tau\oplus j_o)(a)V_n^*-\tau(a)\|=0\,\,\,
{\rm for\,\,\,all}\,\,\, a\in A.
$$
\end{thm}

\begin{proof}
Let $J: B\to l^{\infty}(B)$ be defined by $J(c)=(c,c,...)$ for
$c\in B.$ Define $\tau_{\infty} =J\circ \tau$ and $J_o=J\circ
j_o.$ Let $\omega$ be a free ultrafilter on ${\mathbb N}$ and
$\pi: l^{\infty}(B)\to q_{\omega}(B)$ be the quotient map. It
follows from \ref{IILoinf} that $\pi\circ \tau_{\infty}(A)'$
contains a unital \SCA\, which is isomorphic to ${\cal
O}_{\infty}.$ Denote this \SCA\, by ${\cal O}_{\infty}.$ Let $q\in
{\cal O}_{\infty}$ be a nonzero projection such that $[q]=0$ in
$K_0({\cal O}_{\infty}).$ There is a \SCA\, $C$  of ${\cal
O}_{\infty}$ for which $1_C=q$ and $C\cong {\cal O}_2.$ Put
$\tau_0(a)=q\pi\circ \tau_{\infty}(a).$ So we may view $\tau_0$ is
a unital full \hm\, from $A$ into $\tau_0(A)\otimes {\cal O}_2.$
Since ${\cal O}_2\cong {\cal O}_2\otimes {\cal O}_2$ (by
\cite{R3}), it follows from \ref{IILo22} that $\tau_0\oplus
\pi\circ J_o$ and $\tau_0$ are approximately unitarily equivalent.
Thus $\pi\circ \tau_{\infty}$ and $\pi\circ \tau_{\infty}\oplus
\pi\circ J_o$ are approximately unitarily equivalent.  It follows
from 6.2.5 \cite{Rb} that $\tau$ and $\tau\oplus j_o$ are
approximately unitarily equivalent.

\end{proof}

{\bf Proof of Theorem \ref{MTad}}

\begin{proof}
Since $A$ is separable, there is a unital embedding $j: A \to
{\cal O}_2,$ by 2.8  of \cite{KP}. Since $B$ has property (P2),
there is a full monomorphism $\sigma: {\cal O}_2\to B.$ Define
${\bar j}=\sigma\circ j.$
 Note ${\bar j}$ is full.
Let $\ep>0$ and ${\cal F}\subset A$ be a finite subset. It follows
from Theorem 3.9 of \cite{Lnuct}
 that there is an integer $n$ and a
unitary $v\in M_{n+1}(B)$ such that
\begin{eqnarray}
\|v^*{\rm diag}(h_1(a), {\bar j}(a),{\bar j}(a),...,{\bar
j}(a))v-{\rm diag}(h_2(a), {\bar j}(a),{\bar j}(a),...,{\bar
j}(a))\|<\ep/4
\end{eqnarray}
for all $a\in {\cal F}.$ On the other hand, by
\ref{IIPo2}, there is an isometry $u\in M_n(\pi\circ \sigma({\cal
O}_2))$ with $uu^*=1_{M(C)/C}$ such that
\begin{eqnarray}
\|u^*{\bar j}(a) u-{\rm diag}({\bar j}(a),{\bar j}(a),...,{\bar j}(a))\|<\ep/4
\end{eqnarray}
for $a\in {\cal F}.$  Thus, we obtain an isometry $w\in M_2(B)$
with $ww^*=1_{B}$ such that
\begin{eqnarray}
\|w^*{\rm diag}(h_1(a), j(a))w-{\rm diag}(h_2,
j(a))\|<\ep/2\,\,\,\,{\rm for\,\,\,all}\,\,\, a\in {\cal F}.
\end{eqnarray}
By applying \ref{IIITab}, we obtain a partial isometry $z\in B$ such that
$z^*h_1(1_A)z=h_2(1_A),$ $zh_2(1_A)z^*=h_1(1_A)$ and
\begin{eqnarray}
\|z^*h_1(a)z-h_2(a)\|<\ep\,\,\,\,{\rm for\,\,\, all}\,\,\,a \in
{\cal F}.
\end{eqnarray}
\end{proof}

\begin{Remark}
{\rm
If both $h_1$ and $h_2$ are unital, it is clear that
$z$ can be chosen to be unitary. If one of them is unital and the other
is not, $z$ can never be unitary.
Suppose that both are not unital. Since $B$ has property (P1),(P2) and (P3),
we obtain full ${\cal O}_2$ embeddings into $h_1(1_A)Bh_1(1_A)$
and $h_2(1_A)Bh_2(1_A).$ Therefore there is a projection $e\le h_1(1_A)$
such that $h_1(1_A)$ is equivalent to $h_1(1_A)-e$ and $e$ is a full
projection. So there is a partial isometry $v\in B$ such that
$v^*v=h_1(1_A)$ and $vv^*=h_1(1_A)-e.$
Thus $1-{\rm ad}\, v^*\circ h_1(1_A)$ is full. Similarly, there is
a partial isometry $w\in B$ with $w^*w=h_2(1_A)$ such that
$1-{\rm ad}\, w^*\circ h_2(1_A)$ is full.
Now apply \ref{MTad} to the case that $A={\mathbb C}.$
we know that $1-{\rm ad}v^*\circ h_1(1_A)$ and
$1-{\rm ad}w^*\circ h_2(1_A)$ are equivalent. This implies
that we can choose $z$ to be unitary in the proof of \ref{MTad}.
}

\end{Remark}

\begin{Cor}\label{CPI}
Theorem \ref{MTad} also holds for the case that 
$B=q_{\infty}(\{C_n\}),$ where
each $C_n$ is a unital purely infinite simple \CA s. 
\end{Cor}

\begin{proof}
It is clear that $B$ has property (P1) and (P2). From the proof of
\ref{MTad} above, we only need an absorbing lemma \ref{IIITab} for
this $B.$ Let $\tau: A\to B$ be a full monomorphism and $j_0: A\to
{\cal O}_2\to B$ be a full embedding of $A$ into $B$ which factors
through ${\cal O}_2.$ So we may write $j_0=\Phi\circ j,$ where $j:
A\to {\cal O}_2$ is a monomorphism and $\Phi: {\cal O}_2\to B$ is
a full \hm.  Let $L: A\to l^{\infty}(\{C_n\}$ be a \morp\, for
which $\pi\circ L=\tau,$ where $\pi: l^{\infty}(\{C_n\})\to
q_{\infty}(\{C_n\})$ is the quotient map. Write $L=\{L_n\},$ where
$L_n: A\to C_n$ is a \morp. Let $\phi_n: {\cal O}_2\to C$ such
that $\pi\circ \{\psi_n\}=\Phi.$ Denote by $D_n$ the separable
unital purely infinite simple \CA\, containing $L_n(A)$ and
$\psi_n({\cal O}_2).$ Then $q_{\infty}(\{D_n\})\subset B$ and
$\tau: A\to q_{\infty}(\{C_n\})$ and $j_0: A\to {\cal O}_2\to
q_{\infty}(\{C_n\}).$ Thus one applies 7.5 of \cite{Lnpro}.
\end{proof}

{\bf Proof of Proposition \ref{MPg1}}

\begin{proof}
Let $h_1: A\to B\otimes {\cal K}$ be a \hm. It follows from 4.5 in
\cite{Lnpro} that there is a sequence of asymptotically
multiplicative \morp s $\{\phi_n\}$ from $A$ to $B\otimes {\cal
K}$ and a sequence of unitaries $u_n\in {\widetilde{B\otimes {\cal
K}}}$ such that
\begin{eqnarray}
\lim_{n\to\infty}\|(h\oplus \phi_n)(a)-{\rm ad}\, u_n\circ
j(a)\|=0\rforal a\in A.
\end{eqnarray}
Since $B$ has property (P2), it is easy to see that we may assume
that $\phi_n$ maps $A$ into $B$ and $u_n$ are unitaries in $B.$ It
follows from 6.5 in \cite{Lnsemi} that, for each $k,$ there exists
a sequence of unitaries $v_n(k)\in M_2(B)$ such that
\begin{eqnarray}
\lim_{n\to\infty}\|v_{n}(k)^*(\phi_n(a)\oplus j_o(a) )
v_n(k)-(\phi_{n+k}(a)\oplus j_o(a))\|=0\rforal a\in A.
\end{eqnarray}
It follows from 4.7 of \cite{Lnpro} that there exists a \hm\,
$h_1: A\to M_2(B)$ and a sequence of unitaries $w_n\in M_2(B)$
such that
\begin{eqnarray}
\lim_{n\to\infty}\|{\rm ad}\, w_n\circ h_1(a)-(\phi_n(a)\oplus
j_o(a))\|=0\rforal a\in A.
\end{eqnarray}
By applying the fact that $B$ has property (P2) and applying
\ref{IIPo2}, we obtain a sequence of isometries  $z_n\in M_3(B)$
with $z_nz_n^*=j_o(1_A)$ such that
\begin{eqnarray}
\lim_{n\to\infty}\|(h\oplus h_1\oplus
j_o)(a)-z_n^*j_o(a)z_n\|=0\rforal a\in A.
\end{eqnarray}
It follows that $[h_1]=-[h]$ in $H(A,B).$
\end{proof}

{\bf Proof of \ref{MPg2}}

\begin{proof}
The corollary follows immediately from \ref{MPg1} and
\ref{IIITab}.
\end{proof}

\section{Classification of full extensions}

\begin{Def}\label{DKunder}
{\rm Let $C_n$ be a commutative \CA\, with $K_0(C_n)={\mathbb
Z}/n{\mathbb Z}$ and $K_1(C_n)=0.$ Suppose that $A$ is a \CA. Put
$K_i(A, {\mathbb  Z}/k{\mathbb  Z})=K_i(A\otimes C_k)$ (see
\cite{Sch1}).  One has the following six-term exact sequence (see
\cite{Sch1}):
$$
\begin{array}{ccccc}
K_0(A) & \to & K_0(A, {\mathbb Z}/k{\mathbb Z}) & \to & K_1(A)\\
\uparrow_{{\bf k}} & & & & \downarrow_{{\bf k}}\\
K_0(A) & \leftarrow & K_1(A,{\mathbb Z}/k{\mathbb Z}) & \leftarrow
&
K_1(A)\,.\\
\end{array}
$$
In \cite{DL}, $K_i(A, {\mathbb Z}/n{\mathbb Z})$ is identified
with $KK^i({\mathbb I}_n, A)$ for $i=0,1$.
 As in \cite{DL}, we use the notation
$$
{\underline K}(A)=\oplus_{i=0,1, n\in {\mathbb  Z}_+}
K_i(A;{\mathbb Z}/n{\mathbb Z}).
$$
By ${ \rm Hom}_{\Lambda}({\underline K}(A),{\underline K}(B))$ we
mean all \hm s from ${\underline K}(A)$ to ${\underline K}(B)$
which respect the direct sum decomposition and the so-called
Bockstein operations (see \cite{DL}). It follows from the
definition in \cite{DL} that if $x\in KK(A,B),$ then the Kasparov
product $KK^{i}({\mathbb I}_n, A)\times x$ gives an element in
$KK^i({\mathbb I}_n, B)$ which we identify with ${ \rm
Hom}(K_i(A,{\mathbb Z}/n{\mathbb Z}), K_0(B,{\mathbb Z}/n{\mathbb
Z})).$ Thus one obtains a map $\Gamma: KK(A,B)\to { \rm
Hom}_{\Lambda}(\underline{K}(A),\underline{K}(B)).$ It is shown by
Dadarlat and Loring (\cite{DL}) that if $A$ is in ${\cal N}$
then, for any $\sigma$-unital \CA\, $B,$ the map $\Gamma$ is
surjective and ${\rm ker}\,\Gamma=Pext(K_*(A),K_*(B)).$ In
particular,
$$
\Gamma: KL(A,B)\to{ \rm Hom}_{\Lambda}({\underline
K}(A),{\underline K}(B))
$$
is an isomorphism.
It is shown in \cite{Lnuct} that if $A$ satisfies
AUCT,  then $\Gamma$ is also an isomorphism from $KL(A,B)$
onto ${ \rm Hom}_{\Lambda}({\underline
K}(A),{\underline K}(B)).$

}

\end{Def}

\begin{Lem}\label{IIL2p}
Let $B$ be a unital \CA\, which admits a full ${\cal O}_2$
embedding and let $G_i$ be
a countable subgroup of $K_i(B)$ $(i=0,1).$ There exists a unital
separable \CA\, $B_0\subset B$ which has a full ${\cal O}_2$ embedding
such that
$K_i(B_0)\supset G_i$ and $j_{*i}=\id_{K_i(B_0)},$ where $j:
B_0\to B$ is the embedding.
\end{Lem}

\begin{proof}
Let $p_1,...,p_n,...$ be projections and $u_1,u_2,...,u_n,...$ be
unitaries in  $\cup_{k=1}^{\infty} M_k(B)$ such that $\{p_n\}$ and
$\{u_n\}$  generates of $G_0$ and $G_1,$ respectively. There is a
countable set $S$ such that
$$
p_n, u_n \in \cup_{n=1}^{\infty}\{(a_{ij})_{n\times n}\in M_n(B): a_{ij}\in S\}
$$
Let $j_o: {\cal O}_2\to B$ be a full embedding. Let $p=j(1_{{\cal
O}_2})$ and $x_1,x_2,...,x_m\in B$ such that $\sum_{i=1}^m
x_i^*px_i=1.$ Let $B_1$ be the  unital separable \SCA\, generated
by $S$, $\{x_1,x_2,...,x_m\}$ and $j({\cal O}_2).$  Then $B_1$ has
a full ${\cal O}_2$ embedding and $p_n, u_n\in \cup_{k=1}^{\infty}
M_k(B_1)$ for all $n.$ Note that $K_i(B_1)$ is countable. The
embedding $j_1: B_1\to B$ gives \hm s $(j_1)_{*i}: K_i(B_1)\to
K_i(B).$ Let $F_{1,i}$ be the subgroup of $K_0(B_1)$ generated by
$\{p_n\}$ and $\{u_n\},$ respectively. It is clear that
$(j_1)_{*i}$ is injective on $F_{1,i},$ $i=0,1.$  In particular,
the image of $(j_1)_{*i}$ contains $G_i,$ $i=0,1.$ Let
$N_{1,i}'={\rm ker}(j_1)_{*i}$ and let $N_{1,i}$ be the set of all
projections (if $i=0$), or unitaries (if $i=1$) in
$\cup_{k=1}^{\infty}M_k(B_1)$ which have images in $N_{1, i}'.$
 Let $\{p_{1,n}\}$ be a dense subset of projections in
 $\cup_{k=1}^{\infty}M_k(B_1).$ There are countable pairs
 of projections $\{e_n, e_n'\}$ in $\{p_{1,n}\}$ such that
 $[e_n]=[e_n']$ in $K_0(B).$
 There are $w_n\in \cup_{k=1}^{\infty} M_k(B)$ such that
 $w_n^*w_n=e_n\oplus 1_{k(n)}$ and $w_nw_n^*=e_n'\oplus 1_{k(n)}.$

 Let
 $\{u_{1,n}\}$ be a dense subset of unitaries in
 $\cup_{k=1}^{\infty} M_k(B_1).$
For each $u_{1,n},$ there are unitaries $z_{1,n,k}\in
\cup_{j=1}^{\infty}M_j(B),$ $ k=1,2,...,m(n)$ such that
$$
\|z_{1,n,1}-1\|<1/2, \|z_{1,n,m(n)}-u_{1,n}\|<1/2 \andeqn
\|z_{1,n,k}-z_{1,n,k+1}\|<1/2,
$$
$k=1,2,...,m(n),$ $n=1,2,....$ Let $B_2$ be a separable unital
\CA\, containing $B_1$ such that $\cup_{k=1}^{\infty} M_k(B_2)$
contains all $\{w_{1,n}\}$ and $\{z_{1,n,k}\}.$ Note that there is
a full embedding of ${\cal O}_2$ to $B_2.$
 Note also that if $p, q\in \cup_{k=1}^{\infty}M_k(B_1)$ are
projections so that $[p]-[q]\in N_{1,0}$ then $[p]-[q]=0$ in
$K_0(B_2).$ Similarly, if $u\in B_1$ and $[u]\in N_{1,1},$ then
$[u]=0$ in $B_2.$ Suppose that $B_l$ has been  constructed. Let
$j_l: B_l\to B$ be the embedding. Let $N_{l, i}={\rm
ker}(j_l)_{*i},$ $i=0,1.$ As before, we obtain a unital separable
\CA\, $B_{l+1}\supset B_l$ such that every pair projections
$p,\,q\in \cup_{k=1}^{\infty}M_k(B_l)$ with $[p]-[q]\in N_{l,0}$
has the property that $[p]=[q]$ in $K_0(B_{l+1}),$ and every
unitary $u\in B_l$ with $[u]\in N_{l,1}$ has the property that
$[u]=0$ in $K_1(B_{l+1}).$ Let $B_0$ be the closure of
$\bigcup_{l=1}^{\infty} B_l.$ Note $B_0$ admits a full ${\cal
O}_2$ embedding. Note  also that $B_0$ is separable. Let $j:
B_0\to B$ be the embedding.

We claim that $j_{*i}$ is injective. Suppose that $p\,,q\in
M_k(B_0)$ is a pair of projections for which  $[p]-[q]\in {\rm
ker}j_{*0}$ and $[p]-[q]\not=0$ in $B_0.$ Without loss of
generality, we may assume that $p\,,q\in M_k(B_l)$ for some large
integer $l.$ Then $[p]-[q]$ must be in the ${\rm ker}(j_l)_{*0}.$
By the construction, $[p]-[q]=0$ in $K_0(B_{l+1}).$ This would
imply that $[p]-[q]=0$ in $K_0(B_0).$ Thus $j_{*0}$ is injective.
An exactly same argument shows that $j_{*1}$ is also injective.
The lemma then follows.
\end{proof}

\begin{Lem}\label{IILsix}
Let $B$ be a unital \CA\, which admits a full ${\cal O}_2$
embedding. Suppose that $G_i\subset K_i(B)$ and $F_i(k)\subset
K_i(B,{\mathbb Z}/k{\mathbb Z})$ are countable subgroups such that
the image of $F_i(k)$ in
$K_{i-1}(B)$ is contained in $G_{i-1}$ $(i=0,1, k=2,3,...)$. Then
there exists a separable unital \CA\, $C\subset B$ which admits a
full ${\cal O}_2$ embedding such that $K_i(C)\supset G_i,$
$K_i(C,{\mathbb Z}/k{\mathbb Z}) \supset F_i(k)$ and the embedding
$j: C\to B$ induces an injective map $j_{*i}: K_i(C)\to K_i(B)$
and an injective map $j_*: K_i(C, {\mathbb Z}/k{\mathbb Z})\to
K_i(B,{\mathbb Z}/k{\mathbb Z}),$
 $k=2,3,....$
\end{Lem}

\begin{proof}
It follows from  \ref{IIL2p} that there is a separable unital
\CA\, $C_1$ which admits a full ${\cal O}_2$ embedding such that
$K_0(C_1)\supset G_0$ and $K_1(C_1)\supset G_1$ and $j$ induces an
identity map on $K_0(C_1)$ and $K_1(C_1),$ where $j: C_1\to B$ is
the embedding. Fix $k,$ and let $\{x\in K_i(C_1):
kx=0\}=\{g_1^{(i)},g_2^{(i)},...,\}.$ Suppose that
$\{s_1^{(i)},s_2^{(i)},...,\}$ is a subset of $K_{i-1}(B,{\mathbb
Z}/k{\mathbb Z}) $ such that the map from $K_{i-1}(B,{\mathbb
Z}/k{\mathbb Z})$ to $K_{i}(B)$ maps $s_j^{(i)}$ to $g_j^{(i)}.$
For each $z^{(i)}\in K_{i-1}(C_1, {\mathbb Z}/k{\mathbb Z}),$
there is $s_j^{(i)}$ such that $z^{(i)}-s_j^{(i)}\in
K_i(B)/kK_i(B).$ Since $K_i(C_1)$ is countable, the set of all
possible $z^{(i)}-s_j^{(i)}$ is countable. Thus one obtains a
countable subgroup $G_i'$ which contains $K_i(C_1)$ for which
$G_i'/kK_i(B)$ contains the above the mentioned countable set as
well as $F_i(k)\cap ( K_i(B)/kK_i(B))$ for each $k.$ Since
countably many countable sets is still countable, we obtain a
countable subgroup $G_i^{(2)}\subset K_i(B)$ such that $G_i^{(2)}$
contains $G_i'$ and $kK_i(B)\cap G_i^{(2)}= kG_i^{(2)},$
$k=1,2,...,$ and $i=0,1.$ Note also $F_i(k)\cap
(K_i(B)/kK_i(B))\subset G_i^{(2)}/kK_i(B).$ By applying
\ref{IIL2p}, we obtain a separable unital \CA\, $C_2\supset C_1$
such that $K_i(C_2)\supset G_i^{(2)}$ and an embedding from $C_2$
to $B$ gives an injective map on $K_i(C_2),$ $i=0,1.$ Repeating
what we have done above, we obtain an increasing sequence of
countable subgroups $G_i^{(n)}\subset K_i(B)$ such that
$G_i^{(n)}\cap kK_i(B)=kG_i^{(n)}$ for all $k$ and $i=0,1$ and an
increasing sequence of separable  \SCA s $C_n$ such that
$K_i(C_n)\supset G_i^{(n)}$ and embeddings from $C_n$ into $B$
giving injective maps on $K_i(C_n),$ $i=0,1,$ and $n=1,2,....$
Moreover $F_i^{(k)}\cap (K_i(B)/kK_i(B))\subset K_i(C_n)/kK_i(B).$
Let $C$ denote the closure of $\bigcup_nC_n$ and $j: C\to B$ be
the embedding. Then $C$ is a separable  unital \CA\, and $j_{*i}$
is an injective map, $i=0,1.$ Since $C\supset C_1$ and $C_1$ is
unital, $C$ admits a full ${\cal O}_2$ embedding. We claim that
$K_i(C)\cap kK_i(B)=kK_i(C),$ $k=1,2,...,$ and $i=0,1.$ Note that
$K_i(C)=\cup_nG_i^{(n)}.$ Since $G_i^{(n)}\cap kK_i(B)=kG_i^{(n)}
\subset kK_i(C),$ we see that $K_i(C)\cap kK_i(B)=kK_i(C),$
$i=0,1.$ Thus $K_i(C)/kK_i(C)=K_i(C)/kK_i(B).$ Since
$K_i(C)/kK_i(B)\supset F_i^{(k)}\cap(K_i(B)/kK_0(B)),$ we conclude
also that $K_i(C,{\mathbb Z}/k{\mathbb Z})$ contains $F_i(k).$
Since $j_{*i}$ is injective, $j$ induces an injective map  from
$K_i(C)/kK_i(C)$ into $K_i(B)/kK_i(B)$ for all integer $k\ge 1.$
Using this fact  and the fact that $j_{*i}: K_i(C)\to K_i(B)$ is
 injective,
  by chasing the following commutative diagram
 $$
{\small \put(-160,0){$K_0(C)$} \put(0,0){$K_0(C,{\bf Z}/k{\bf
Z})$} \put(180,0){$K_1(C)$} \put(-85,-40){$K_0(B)$}
\put(0,-40){$K_0(B,{\bf Z}/k{\bf Z})$} \put(105,-40){$K_1(B)$}
\put(-85, -70){$K_0(B)$} \put(0,-70){$K_1(B, {\bf Z}/k{\bf Z})$}
\put(105,-70){$K_1(B)$} \put(-160,-110){$K_0(C)$}
\put(0,-110){$K_1(C,{\bf Z}/k{\bf Z})$} \put(180,-110){$K_1(C)$}
\put(-120, 2){\vector(1,0){95}} \put(70,1){\vector(1,0){95}}
\put(-123,-3){\vector(1,-1){30}} \put(30,-3){\vector(0,-1){25}}
\put(180,-2){\vector(-1,-1){30}} \put(-45,-38){\vector(1,0){35}}
\put(70,-38){\vector(1,0){25}} \put(-147, -90){\vector(0,1){85}}
\put(-75,-60){\vector(0,1){15}} \put(115, -45){\vector(0,-1){15}}
\put(190,-7){\vector(0,-1){85}} \put(-7,-68){\vector(-1,0){35}}
\put(95,-68){\vector(-1,0){25}} \put(-123,-102){\vector(1,1){30}}
\put(175, -105){\vector(-1,1){30}} \put(30,-100){\vector(0,1){25
}} \put(-5, -108){\vector(-1,0){100}}
\put(170,-108){\vector(-1,0){95}} \put(-112,-12){$j_{*0}$}
\put(15, -15){$j_* $} \put(150, -14){$j_{*1}$}
\put(-130,-92){$j_{*0}$} \put(15, -88){$j_*$} \put(160,
-88){$j_{*1}$} }
$$
one sees that $j$ induces an injective map from $K_i(C, {\mathbb
Z}/k{\mathbb Z}) $ to $K_i(B, {\mathbb Z}/k{\mathbb Z}).$
\end{proof}

\begin{Cor}\label{IIfu}
Without assuming that $B$ has a full ${\cal O}_2$ embedding, both
\ref{IILsix} and \ref{IIL2p} hold if we do not require that $C$
(or $B_0$) has a full ${\cal O}_2$ embedding.
\end{Cor}

{\bf Proof of Theorem \ref{MT2}}

\begin{proof}
By \ref{MT1}, it suffices to show that, for each $x\in KL(A,
M(B)/B),$ there is a full monomorphism $h: A\to M(B)/B$ such that
$[h]=x.$ Put $Q=M(B)/B.$ Since $A$ satisfies the AUCT, we may view
$x$ as an element in ${\rm Hom}_{\Lambda}(\underline{K}(A),
\underline{K}(Q)).$ Note that $K_i(A)$ is a countable abelian
group ($i=0,1$). Let $G_{0}^{(i)}=\gamma(x)(K_i(A)),$ $i=0,1,$
where $\gamma: {\rm Hom}_{\Lambda}(\underline{K}(A),
\underline{K}(Q)) \to {\rm Hom}(K_*(A),K_*(Q))$ is the surjective
map. Then $G_{0}^{(i)}$ is a countable subgroup of $K_i(Q),$
$i=0,1.$ Consider the following commutative diagram:
$$
{\small \put(-160,0){$K_0(A)$} \put(0,0){$K_0(A,{\bf Z}/k{\bf
Z})$} \put(180,0){$K_1(A)$} \put(-85,-40){$K_0(Q)$}
\put(0,-40){$K_0(Q,{\bf Z}/k{\bf Z})$} \put(105,-40){$K_1(Q)$}
\put(-85, -70){$K_0(Q)$} \put(0,-70){$K_1(Q, {\bf Z}/k{\bf Z})$}
\put(105,-70){$K_1(Q)$} \put(-160,-110){$K_0(A)$}
\put(0,-110){$K_1(A,{\bf Z}/k{\bf Z})$} \put(180,-110){$K_1(A)$}.
\put(-120, 2){\vector(1,0){95}} \put(70,1){\vector(1,0){95}}
\put(-123,-3){\vector(1,-1){30}} \put(30,-3){\vector(0,-1){25}}
\put(180,-2){\vector(-1,-1){30}} \put(-45,-38){\vector(1,0){35}}
\put(70,-38){\vector(1,0){25}} \put(-147, -90){\vector(0,1){85}}
\put(-75,-60){\vector(0,1){15}} \put(115, -45){\vector(0,-1){15}}
\put(190,-7){\vector(0,-1){85}} \put(-7,-68){\vector(-1,0){35}}
\put(95,-68){\vector(-1,0){25}} \put(-123,-102){\vector(1,1){30}}
\put(175, -105){\vector(-1,1){30}} \put(30,-100){\vector(0,1){25}}
\put(-5, -108){\vector(-1,0){100}}
\put(170,-108){\vector(-1,0){95}} \put(-111,-14){$\gamma(x)$}
\put(15, -15){$\times x $} \put(148, -14){$\gamma(x)$}
\put(-130,-90){$\gamma(x)$} \put(15, -88){$\times x$} \put(160,
-88){$\gamma(x)$} }
$$
It follows from \ref{IILsix} that there is a unital
\CA\, $C\subset Q$ which has a full ${\cal O}_2$ embedding
such that $K_i(C)\subset G_0^{(i)},$
$K_i(C)\cap kK_i(Q)=kK_i(C),$ $k=1,2,...,$ and $i=0,1,$ and
the embedding $j: C\to Q$ induces injective maps on $K_i(C)$ as well
as on $K_i(C,{\mathbb Z}/k{\mathbb Z})$ for all $k$ and $i=0,1.$
Moreover $K_i(C,{\mathbb Z}/k{\mathbb Z})\supset
(\times x)(K_i(A,{\mathbb Z}/k{\mathbb Z}))$ for $k=1,2,...$ and $i=0,1.$
We
have the following commutative diagram:
$$
{\small \put(-160,0){$K_0(A)$} \put(0,0){$K_0(A,{\bf Z}/k{\bf
Z})$} \put(180,0){$K_1(A)$} \put(-85,-40){$K_0(C)$}
\put(0,-40){$K_0(C,{\bf Z}/k{\bf Z})$} \put(105,-40){$K_1(C)$}
\put(-85, -70){$K_0(C)$} \put(0,-70){$K_1(C, {\bf Z}/k{\bf Z})$}
\put(105,-70){$K_1(C)$} \put(-160,-110){$K_0(A)$}
\put(0,-110){$K_1(A,{\bf Z}/k{\bf Z})$} \put(180,-110){$K_1(A)$}
\put(-120, 2){\vector(1,0){95}} \put(70,1){\vector(1,0){95}}
\put(-123,-3){\vector(1,-1){30}}
\put(180,-2){\vector(-1,-1){30}} \put(-45,-38){\vector(1,0){35}}
\put(70,-38){\vector(1,0){25}} \put(-147, -90){\vector(0,1){85}}
\put(-75,-60){\vector(0,1){15}} \put(115, -45){\vector(0,-1){15}}
\put(190,-7){\vector(0,-1){85}} \put(-7,-68){\vector(-1,0){35}}
\put(95,-68){\vector(-1,0){25}} \put(-123,-102){\vector(1,1){30}}
\put(175, -105){\vector(-1,1){30}}
\put(-5, -108){\vector(-1,0){100}}
\put(170,-108){\vector(-1,0){95}} \put(-112,-12){$\gamma(x)$}
 \put(145, -12){$\gamma(x)$}
\put(-135,-92){$\gamma(x)$}  \put(160, -88){$\gamma(x)$} }
$$
We will add two more maps on the above diagram.
 From
 the fact that  the image of $K_i(A, {\mathbb Z}/k{\mathbb Z})$
 under $\times x$ is contained in $K_i(C,{\mathbb Z}/k{\mathbb
 Z}),$ ($k=2,3,...,$ $i=0,1$),
  we obtain two maps $\beta_i:K_i(A, {\mathbb Z}/k{\mathbb
 Z})\to K_i(C,{\mathbb Z}/k{\mathbb
 Z}),$ $k=2,3,...,$ $i=0,1$ such that $j_*\circ \beta_i=\times x$
 and obtain the
 following commutative diagram:
$$
{\small \put(-160,0){$K_0(A)$} \put(0,0){$K_0(A,{\bf Z}/k{\bf
Z})$} \put(180,0){$K_1(A)$} \put(-85,-40){$K_0(Q)$}
\put(0,-40){$K_0(C,{\bf Z}/k{\bf Z})$} \put(105,-40){$K_1(C)$}
\put(-85, -70){$K_0(C)$} \put(0,-70){$K_1(C, {\bf Z}/k{\bf Z})$}
\put(105,-70){$K_1(C)$} \put(-160,-110){$K_0(A)$}
\put(0,-110){$K_1(A,{\bf Z}/k{\bf Z})$} \put(180,-110){$K_1(A)$}
\put(-120, 2){\vector(1,0){95}} \put(70,1){\vector(1,0){95}}
\put(-123,-3){\vector(1,-1){30}} \put(30,-3){\vector(0,-1){25}}
\put(180,-2){\vector(-1,-1){30}} \put(-45,-38){\vector(1,0){35}}
\put(70,-38){\vector(1,0){25}} \put(-147, -90){\vector(0,1){85}}
\put(-75,-60){\vector(0,1){15}} \put(115, -45){\vector(0,-1){15}}
\put(190,-7){\vector(0,-1){85}} \put(-7,-68){\vector(-1,0){35}}
\put(95,-68){\vector(-1,0){25}} \put(-123,-102){\vector(1,1){30}}
\put(175, -105){\vector(-1,1){30}} \put(30,-100
){\vector(0,1){25}} \put(-5, -108){\vector(-1,0){100}}
\put(170,-108){\vector(-1,0){95}} \put(-112,-12){$\gamma(x)$}
\put(15, -15){$\beta_0 $} \put(145, -12){$\gamma(x)$}
\put(-135,-92){$\gamma(x)$} \put(15, -88){$\beta_1$} \put(160,
-88){$\gamma(x)$} }
$$
Consider the following commutative diagram:
$$
\begin{array}{ccccccc}
\to & K_i(A,{\mathbb Z}/mn{\mathbb Z}) & \to & K_i(A, {\mathbb Z}/n{\mathbb Z})&
\to & K_{i-1}(A, {\mathbb Z}/m{\mathbb Z}) &\to \\
& \downarrow & & \downarrow && \downarrow\\
\to & K_i(Q,{\mathbb Z}/mn{\mathbb Z}) & \to & K_i(Q, {\mathbb Z}/n{\mathbb Z})&
\to & K_{i-1}(Q, {\mathbb Z}/m{\mathbb Z}) &\to \cr
\end{array}
$$
Since $j_*\circ \beta_i=\times x$ and
all vertical maps in the following diagram is injective
$$
\begin{array}{ccccccc}
\to & K_i(C,{\mathbb Z}/mn{\mathbb Z}) & \to & K_i(C, {\mathbb Z}/n{\mathbb Z})&
\to & K_{i-1}(C, {\mathbb Z}/m{\mathbb Z}) &\to \\
& \downarrow & & \downarrow && \downarrow\\
\to & K_i(Q,{\mathbb Z}/mn{\mathbb Z}) & \to & K_i(Q, {\mathbb Z}/n{\mathbb Z})&
\to & K_{i-1}(Q, {\mathbb Z}/m{\mathbb Z}) &\to \cr,
\end{array}
$$
we obtain
the following commutative diagram:
$$
\begin{array}{ccccccc}
\to & K_i(A,{\mathbb Z}/mn{\mathbb Z}) & \to & K_i(A, {\mathbb Z}/n{\mathbb Z})&
\to & K_{i-1}(A, {\mathbb Z}/m{\mathbb Z}) &\to \\
& \downarrow & & \downarrow && \downarrow\\
\to & K_i(C,{\mathbb Z}/mn{\mathbb Z}) & \to & K_i(C, {\mathbb Z}/n{\mathbb Z})&
\to & K_{i-1}(C, {\mathbb Z}/m{\mathbb Z}) &\to \cr
\end{array}
$$
Thus we obtain an element $y\in KL(A, C)$ such that $y\times
[j]=x.$ Since $A$ satisfies the AUCT, one checks that $KL(A,
C)=KL(A\otimes {\cal O}_{\infty}, C).$ This also follows from the
fact that the unital embedding from $A\to A\otimes {\cal
O}_{\infty}$ gives a $KK$-equivalence (see \cite{Pi}). It follows
from 6.6 and 6.7 in \cite{Lnsemi} that there exists a \hm\, $\phi:
A\otimes {\cal O}_{\infty}\to C\otimes {\cal K}$ such that
$[\phi]=y.$ Define $\psi=\phi|_{A\otimes 1}.$ By the same result
of Pimsner (\cite{Pi}),  one obtains that $[\psi]=y.$ Since $A$ is
unital, we may assume that the image of $\psi$ is in $M_m(C)$ for
some integer $m\ge 1.$  Since $C$ admits a full ${\cal O}_2$
embedding, $C$ has property (P2). Thus $1_m$ is equivalent to a
projection in $C.$ Thus we may further assume that $\psi$ maps $A$
into $C.$
 Put $h_1=j\circ \psi.$ To obtain a full monomorphism, we
note that there is an embedding $\imath: A\to {\cal O}_2$ (see
Theorem 2.8 in \cite{KP}). Since $M(B)/B$ has property (P2), we
obtain a full monomorphism $\psi: {\cal O}_2\to M(B)/B.$ Let
$e=\psi(1_{{\cal O}_ 2}).$ There is a partial isometry $w\in
M_2(M(B)/B)$ such that $w^*w=1_{M(B)/B}$ and $ww^*=1\oplus e.$
Define $h=w^*(h_1\oplus \psi\circ \imath)w.$ One checks that
$[h]=[h_1]=x$ and $h$ is a full monomorphism.
\end{proof}

\begin{Cor}\label{C=}
Let $A$ be a unital separable amenable \CA\, satisfying the AUCT.
Let $B$ be a unital \CA\, which has property (P2). Then, for each
$x\in KL(A,B),$ there is a full monomorphism $h: A\to B$ such that
$[h]=x.$
\end{Cor}

\begin{proof}
In the proof above, we may replace $M(B)/B$ by $B.$
\end{proof}



{\bf Proof of Theorem \ref{MT3}}
\begin{proof}
For the first part of the theorem, it suffices to show that every
essential full extension is absorbing. Let $\tau$ be a such
extension. Following Elliott and Kocerovsky, we will show that
$\tau$ is purely large. Denote $E=\tau^{-1}(A).$ Choose $c\in
E\setminus C.$ Then, by \ref{Ifull2},  $c$ is a full element.
Since $M(C)$ has property (P1), there exists $x\in M(B)$ such that
$x^*cc^*x=1.$ Therefore there exists a projection $p\le cc^*$ for
which there is $v\in M(B)$ such that $v^*v=1$ and $vv^*=p.$ Note
$\overline{cBc^*}= \overline{cM(B)c^*}\cap B.$ So $pBp\subset
\overline{cBc^*}.$ Now $v^*pBpv=B.$ So $pBp$ is stable and $pBp$
is full. Thus $\tau$ is purely large. So it is absorbing. The last
part of the theorem follows from the next corollary.
\end{proof}

\begin{Cor}\label{IVcst1}
Let $A$ be a separable unital amenable \CA, $C$ be a unital
\CA\, and $B=C\otimes {\cal K}.$
 Then $Ext(A,B)$ is the same set as unitary equivalence
classes of essential full extensions of $A$ by $B.$
\end{Cor}

\begin{proof}
It suffices to show that given any element $x\in Ext(A,B),$ there
exists an essential full extension $\tau: A\to M(B)/B$ so that
$[\tau]=x.$ There exists a $\tau_1: A\to M(B)/B$ such that
$[\tau_1]=x.$ Take a monomorphism $j: A\to {\cal O}_2$ (see
\cite{KP}). Let $h: {\cal O}_2\to M({\cal K})$ be a monomorphism
(given by a faithful representation of ${\cal O}_2$ on a separable
Hilbert space). Let $\phi: M({\cal K})\to M(B)$ be the standard
unital embedding and $\pi: M(B)\to M(B)/B$ be the quotient map.
Then $\tau_2=\pi\circ \phi\circ h\circ j$ gives a full essential
trivial extension. It follows that $\tau=\tau_1\oplus \tau_2$ is
an essential full extension. Since $[\tau_2]=0,$
$[\tau]=[\tau_1]=x.$
\end{proof}

\begin{Remark}\label{Rf}
{\rm Let $B$ be a non-stable, non-unital but $\sigma$-unital \CA.
Suppose that $M(B)/B$ has property (P1), (P2) and (P3), and
suppose that $\tau: A\to M(B)/B$ is an essential full extension.
One should not expect that such extension is purely large in
general. Let $0\to B\to E\to A\to 0$ be an essential full
extension corresponding to $\tau.$ Recall that the extension is
purely large if $cBc^*$ contains a \SCA\, which is stable and
$cBc^*$ is full in $B$ (see \cite{EK}). Given any element $c\in
E\setminus B,$ $\pi(c)$ is full in $M(B)/B.$ But, in general, $c$
need not be full in $M(B),$ nor does$cBc^*$ need to be full in
$B.$ Examples are easily seen in the case that $B=c_0(C),$ where
$C$ is a unital purely infinite simple \CA . Suppose that $0\to
c_0(C)\to E\to A\to 0$ is a full extension and $c'\in E\setminus
c_0(C).$ Write $c'=\{c_n'\}\in l^{\infty}(C).$ Define $c_n=c_n'$
if $n\ge N>1$ and $c_n=0$ if $n\le N.$ Put $c=\{c_n\}.$ Then $c\in
E\setminus c_0(C).$ However, it is clear that $cc_0(C)c^*$ is not
full in $c_0(C).$ By \ref{IIITab}, the full extension $\tau$ is
approximately absorbing in the sense of \ref{IIITab} but not
purely large. It should be also noted that, even if $c^*Bc$ is
full for all $c\in E\setminus B,$ the full extension may not be
purely large. Let $B$ be a nonstable, non-unital but
$\sigma$-unital simple \CA\, with continuous scale
(see\cite{Lncs2} for more examples). Then $B$ may be stably
finite. No hereditary \SCA\, of $B$ contains a stable \SCA. So
none of the essential extensions of a unital separable amenable
\CA\, $A$ by $B$ could be possibly purely large in the sense of
\cite{EK}, nevertheless, all of these extensions are approximately
absorbing in the sense of \ref{IIITab} (and many of them are
actually absorbing; for example, when $A=C(X)$).

}

\end{Remark}


\end{document}